\documentclass[reqno]{tatra10a_sciendo_CC}%%% for authors---format A4 to compile without Logos of Tatra and Siendo
\usepackage{hyperref} % nova verzia hyperrefu 6.78s v Miktex-u 2.7 to predefinovava
\makeatletter \thm@headfont{\bfseries\scshape} \makeatother
\usepackage{latexsym,euscript}
\usepackage{amsbsy,amscd,amsfonts,amsmath,amsopn,amssymb,amstext,amsthm,amsxtra,}
\usepackage{epsf,epsfig,graphics,color,verbatim}
\usepackage{cite, float,multirow}
\usepackage{subfigure}

\usepackage[english,greek]{babel}
\usepackage[T1]{fontenc}
\usepackage[utf8]{inputenc}
\usepackage{textalpha}	
\usepackage{multirow}	
\usepackage{graphicx}
\usepackage{amsmath}
\usepackage{amssymb}
\usepackage{bm}  
\usepackage{dcolumn}
\newcolumntype{d}[1]{D{.}{.}{#1}}

\usepackage{chngcntr}

\newcommand\invisiblesubsection[1]{\refstepcounter{subsection}  \addcontentsline{toc}{subsection}{\protect\numberline{\thesubsection}#1} \sectionmark{#1}}

 \newcommand\startsubsectionnumbering{%
    \makeatletter 
\counterwithin{table}{subsection}
\counterwithin{figure}{subsection}
\renewcommand{\thefigure}{\arabic{section}.\arabic{subsection}.\arabic{figure}}
\renewcommand{\thetable}{\arabic{section}.\arabic{subsection}.\arabic{table}}
    \makeatother}
    
\newcommand\startsectionnumbering{
    \makeatletter
\counterwithin{table}{section}
\counterwithin{figure}{section}
\renewcommand{\thefigure}{\arabic{section}.\arabic{figure}}
\renewcommand{\thetable}{\arabic{section}.\arabic{table}}
    \makeatother}      

\theoremstyle{plain}
\theoremstyle{plain}

\theoremstyle{definition}

\def\Inf{\operatornamewithlimits{inf\vphantom{p}}}
\newcommand{\suma}{\underset{i=1}{\overset{n}\sum}}
\newcommand\scalemath[2]{\scalebox{#1}{\mbox{\ensuremath{\displaystyle #2}}}}
\issueinfo{80}{2021}%%%only for  Tatra
\pagespan{71}{118}%%% number of pages 1...
\received{January 31, 2021}%%%%January 1, 2019
\DOI{10.2478/tmmp-2021-0032}%%%  ??? willbe fulfilled %%%%%%
\begin{document}
\selectlanguage{english}
\title%%%%please, write the first letter only with capital letters, others might be written with small letters.
[What was the river Ister in the time of Strabo?]
{What was the river Ister in the time of Strabo? A mathematical approach}
%%%%%%%%%%

\author[1]{Karol Mikula
        }
\affil[1]{Department of Mathematics and Descriptive Geometry, Faculty of Civil Engineering, Slovak University of Technology, Radlinsk\'{e}ho 11, 810 05 Bratislava, Slovakia}
\affil[2]{Department of Classical and Semitic Philology, Faculty of Arts, Comenius University, Gondova 2, 811 02  Bratislava, Slovakia}
\address[1]{Karol Mikula \\
         Department of Mathematics \\ and Descriptive Geometry\\
         Faculty  of Civil Engineering\\
         Slovak University of Technology\\
         Radlinsk\'eho 11        \\
         810 05--Bratislava \\
         Slovakia}
\email{karol.mikula@stuba.sk  }
%%%%%%%%
\author[1]{Martin Ambroz}
\author[2]{Ren\'ata Moko\v sov\'a}
%\affil{University , City, Country}
\address{Martin Ambroz\\
         Department of Mathematics \\ and Descriptive Geometry\\
         Faculty  of Civil Engineering\\
         Slovak University of Technology\\
         Radlinsk\'eho 11        \\
         810 05--Bratislava \\
         Slovakia}
\email{   martin.ambroz.ml@gmail.com    }
\address{Ren\'ata Moko\v sov\'a\\
         Department of Classical \\and Semitic Philology\\
         Faculty  of Arts\\
         Comenius University\\
         Gondova 2        \\
         811 02--Bratislava \\
         Slovakia}
\email{   renata.mokosova@gmail.com    }
%%%%%%%%%%%%%%%%%%%%%%%%%%%%%%%%%%%%%%%%%%%%%%%%%%%%%%%%%%%%%%%%%%%%%%
%%%\author{...}[1]
%%%%\author[2]
%%%\aithor[3]...
%%%%%%%%%%%%%%%%%%%%%%%%%
\def\shortauthors{MIKULA---AMBROZ---MOKO\v SOV\'A}               %%% if the names are too long and overfullbox appear.
%%%%%%%%%%%%%%%%%%%%%%%%%%%%%%%%%%%%%%%%%%%%%%%%%%%%%%%%%%%%%%%%%%%%%%

\thanks {Supported by the Grants APVV-19-0460 and VEGA 1/0436/20.}

%%%%%%%%%%%%%%%%%%%%%%%%%%%%%%%%%%%%%%%%%%%%%%%%%%%%%%%%%%%%%%%%%%%%%%%%%%%%%%%%%%%%%%%%%%%
\begin{abstract}
In this paper, we introduce a novel method for map registration and apply it to transformation of the river Ister from {\it Strabo's map of the World} to the current map in the World Geodetic System. This transformation leads to the surprising but convincing result that Strabo's river Ister best coincides with the nowadays Tauernbach-Isel-Drava-Danube course and not with the Danube river what is commonly assumed. Such a result is supported by carefully designed mathematical measurements and it resolves all related controversies otherwise appearing in understanding and translation of Strabo's original text. Based on this result we also show that {\it Strabo's Suevi in the Hercynian Forest} corresponds to the Slavic people in the Carpathian-Alpine basin and thus that the compact Slavic settlement was there already at the beginning of the first millennium AD. 
\end{abstract}

%%%%%%%%%%%%%%%%%%%%%%%%%%%%%%%%%%%%%
\maketitle
%%%%%%%%%%%%%%%%%%%%%%%%%%%%%%%%%%%%

\par %%use paragraph \par not \\ \\

\section{Introduction}

\startsectionnumbering
\begin{figure}[h]
 \begin{center}
 \includegraphics[width = 1.\textwidth]{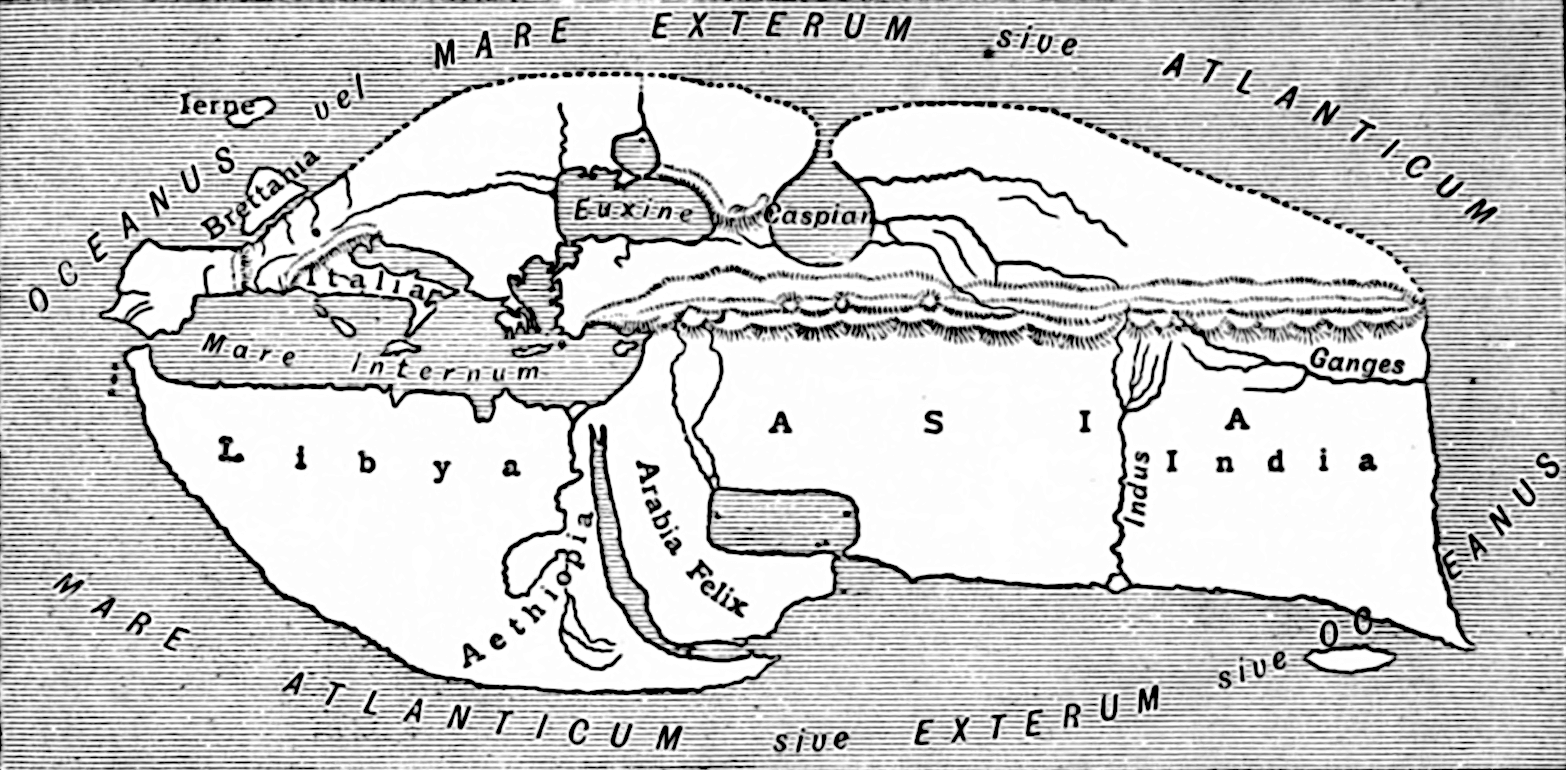}
 \caption{Strabo's map of the World.}
 \label{fig:mapa_strabo}
 \end{center}
\end{figure}

In this paper, we present a novel mathematical model and numerical method for the transformation of geographic maps to each other. To be more precise, we are interested in finding the transformation of a historical map to a current one in the World Geodetic System (WGS) \cite{wgs}. Such a problem is also called map registration. Our mathematical model and numerical method for the map registration is based on two main principles and steps. First, we design and compute locally optimal affine transformations of one map to the other. In this step, every locally optimal affine transformation is found by means of the least square method using a set of clearly identified corresponding points. Then, in the second step, the locally optimal affine transformations are smoothly interpolated/extrapolated to all other points of the map by solving the Laplace equation with suitable boundary conditions. The solution of the Laplace equation is obtained numerically by using the finite difference method on a background grid discretizing the selected rectangular region of interest on the map. After obtaining the final transformation, which we call {\it Locally Affine Globally Laplace} (LAGL) map transformation, we can transform any point of the historical map to the current map and see which places on the current map correspond to geographic objects on the historical map.

The motivation to study this problem mathematically and numerically comes to us by reading Strabo's {\it Geographica} (Στράβωνος Γεωγραφικά)
and studying the related {\it Strabo's map of the World} published in the Encyclopaedia Biblica \cite{EB}, in the section "Geography" on page 1691, see Figure \ref{fig:mapa_strabo}. {\it Strabo's map of the World} is authored by Karl M\"uller, a famous 19th-century historian, classic philologist, geographer, and cartographer, who translated Strabo's {\it Geographica} from Greek to Latin in 1853 \cite{MD}. This book is one of our reading sources because it contains the Greek text without "purposeful" changes found sometimes in other sources and translations. The further very useful Greek text of Strabo's {\it Geographica}, allowing direct translation of Greek words in English are available in Perseus Digital Library \cite{PDL} of Tufts University. It contains the {\it Geographica} edition by August Meineke from 1877 \cite{Meineke} but one has to be careful because at some points it deviates from \cite{MD} and other sources, e.g. by exchanging river names or their transcript with respect to the original. The most recent and very useful source for reading is the English translation of Strabo's work by Duane W. Roller \cite{Roller} where he aims to respect the Strabo's original geographic names and do not translate them to commonly used nowadays terms. Concerning the translation and understanding of the Strabo's work, it is worth to cite Roller's book, Section 5: "there is still the problem of many rare or unique words, extensive paraphrases of earlier authors who are themselves obscure, ambiguities of style and sheer length of the work". There exist some further useful translations which can be found on the internet, e.g. \cite{Lacus}, suitable for an introductory reading. 

Strabo was a Greek geographer who lived from around 63 BC to around AD 24 \cite{Roller}. He lived in Asia Minor, Rome and Alexandria and travelled a lot during his lifetime to collect the information for his work not only by reading preceding sources such as Eratosthenes and others, mainly in the famous Alexandria library. Strabo's {\it Geographica} was first published in 7 BC, collecting all the knowledge from the previous years, and he continues the work until approximately AD 23, during the reigns of the emperors Augustus and Tiberius. Strabo's {\it Geographica} is considered to be one of the rare ancient scientific works in the human history remained to modern times, it is not a historical narrative, but it gives a huge amount of useful quantified information about the known world at the beginning of the first millennium AD. And it is very important to note that we must not apply any later knowledge of the Romans about Europe and the World when reading {\it Geographica}. 
 
As it is announced in the title of the paper, we investigate how the river Ister (Ἴστρος) is transformed from the historical {\it Strabo's map of the World} to the current map of the world.
% in WGS, such as the Google Earth \cite{gEarth} or  the GeoGraphics functionality in Wolfram Mathematica \cite{mathematica}. 
We discovered an astonishing fact that the {\bf Strabo's river Ister}, or better say the river Ister in the times of Strabo, does not correspond to the river Danube on its entire course, but it {\bf perfectly fits with} the nowadays {\bf Tauernbach-Isel-Drava-Danube} course (or if simplified we can just say to the Drava-Danube course). This result is surprising but convincing, supported by carefully designed quantitative mathematical measurements. First, by computing a distance of sources of the current Danube, Drava, Isel and Tauernbach rivers and the source of the transformed Strabo's Ister. Further, by computing the Hausdorff distance of curves representing the respective river courses. And finally, by computing the common length of the compared river courses in prescribed narrow bands. From all these comparisons and from Strabo's writing itself it is clear that the current Tauernbach-Isel-Drava-Danube course gives the best correspondence with the Ister in the Strabo's times. Moreover, this result avoids several, if not all, contradictions which otherwise occur when reading carefully {\it Geographica} and its later translations which used to consider the Ister as the Danube river in the whole its course, i.e. with the sources in the Schwarzwald, Germany. In the sequel, we mention briefly just a few, the most important contradictions, and their solution given by our mathematical result. 

\begin{figure}[ht]
 \begin{center}
 \includegraphics[height = 37mm]{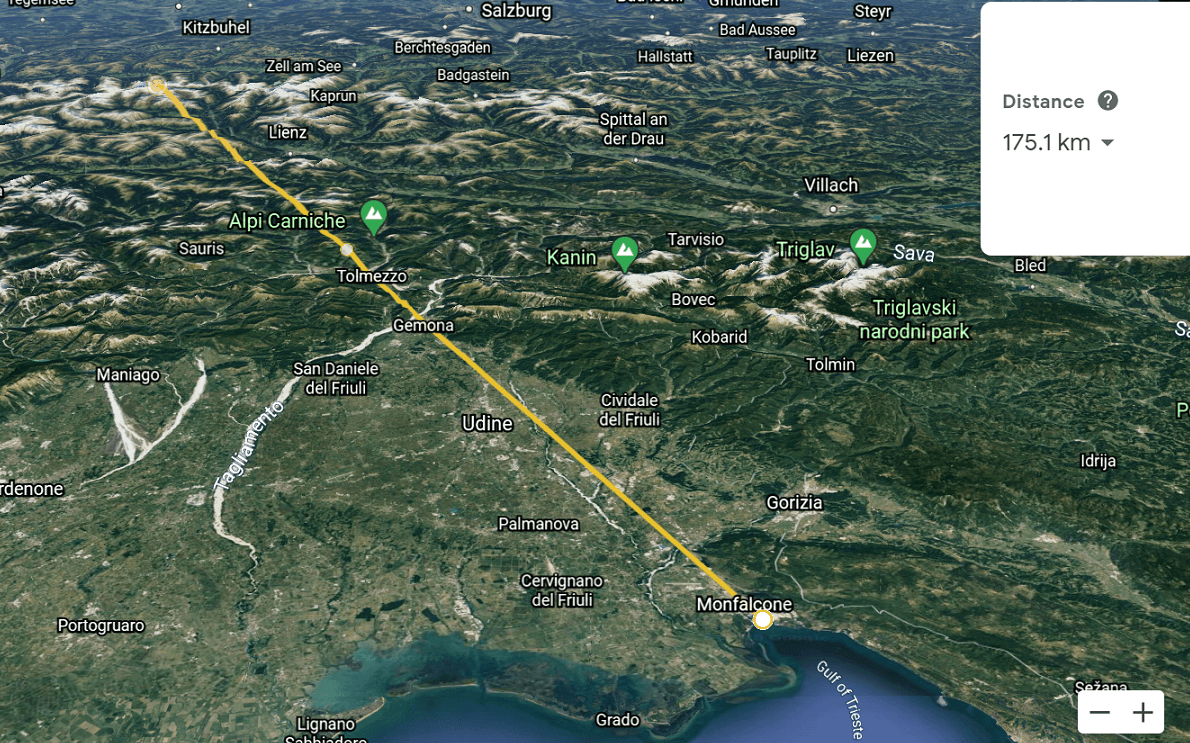}\hskip 2mm
 \includegraphics[height = 37mm]{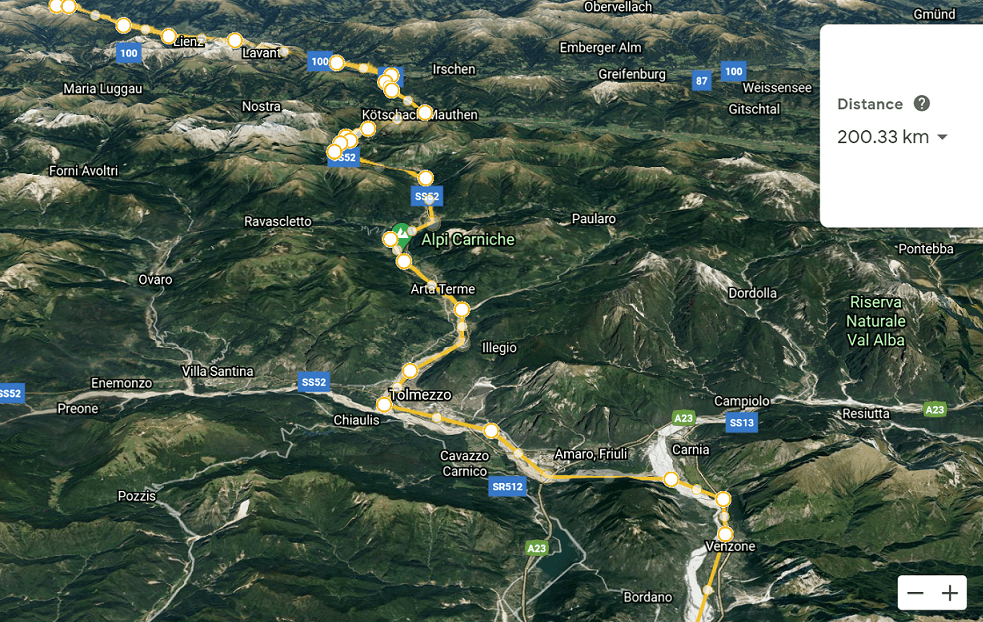}\vskip 2mm
 \includegraphics[height = 56mm]{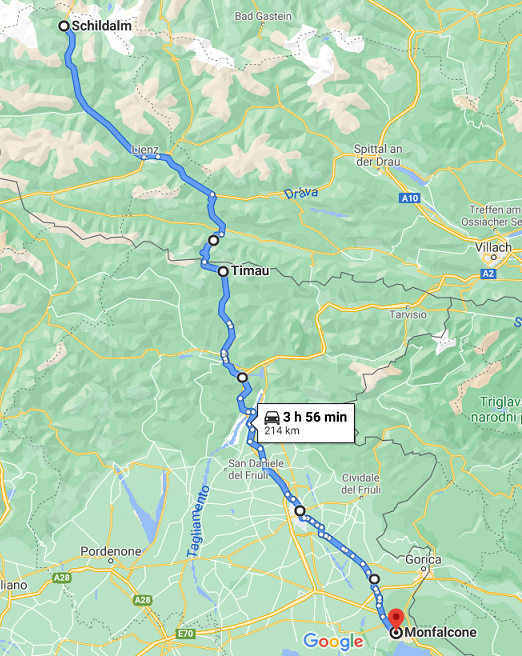}\hskip 6mm
 \includegraphics[height = 56mm]{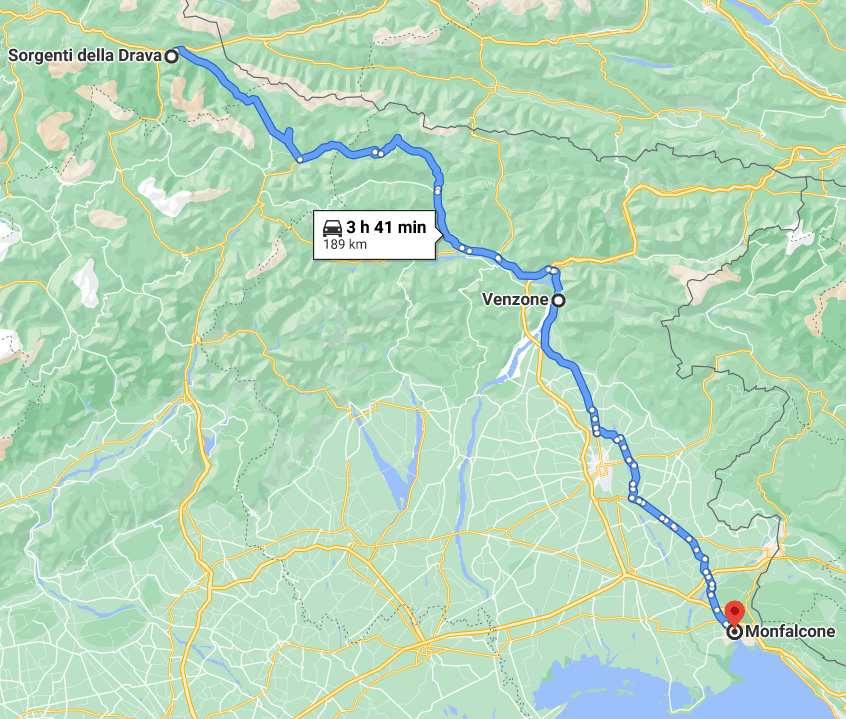}
 \caption{Distances to the recess of Adriatic (Monfalcone). Left up: geodesic line distance from the source of Tauernbach in the High Tauern (175 km). Right up: geodesic distance from the source of the Tauernbach through the passable valleys (200 km). Left down: distance from the source of the Tauernbach by local roads (214 km). Right down: distance from the source of the Drava river by local roads (189 km). Upper images were created by the Google Earth application while bottom images by the Google Maps application.}
 \label{fig:mapa_cesty}
 \end{center}
\end{figure}

First, hardly explainable contradictions are given by the location of the Ister source, its distance  from the recess of Adriatic and by the direction of the Ister course itself, as described by Strabo in the Book 7, Part 1, Section 1 of {\it Geographica}, which we denote by (7-1-1) and other sections are denoted in the same manner. In (7-1-1) Strabo says:

\vskip 3mm\noindent
"Ἴστρος ... ῥέων πρὸς νότον κατ᾽ ἀρχάς, εἶτ᾽ ἐπιστρέφων εὐθὺς ἀπὸ τῆς δύσεως ἐπὶ τὴν ἀνατολὴν καὶ τὸν Πόντον. ἄρχεται μὲν οὖν ἀπὸ τῶν Γερμανικῶν ἄκρων τῶν ἑσπερίων, πλησίον δὲ καὶ τοῦ μυχοῦ τοῦ Ἀδριατικοῦ, διέχων αὐτοῦ περὶ χιλίους σταδίους". 

\begin{figure}[ht]
 \begin{center}
 \includegraphics[width = 1.\textwidth]{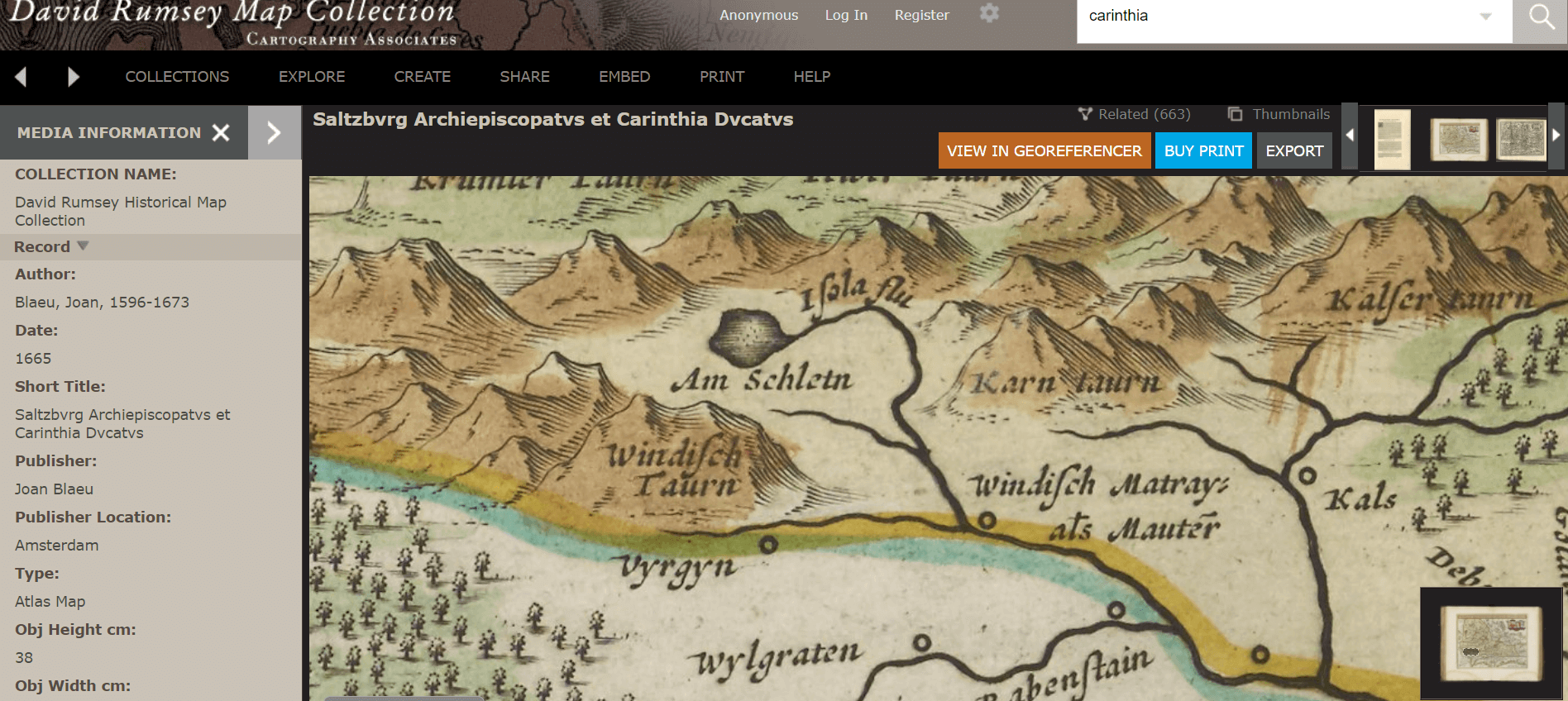}
 \caption{Detail of the map by Joan Blaeu from 1665 \cite{Blaeu} where we see the High Tauern and the Isola flu corresponding to the Tauernbach-Isel stream. The Gro{\ss}venediger corresponds to the Windisch Taurn and the Gro{\ss}glockner to the Kalser Taurn.
}
 \label{fig:isola}
 \end{center}
\end{figure}

\begin{figure}[ht]
 \begin{center}
 \includegraphics[width = 1.\textwidth]{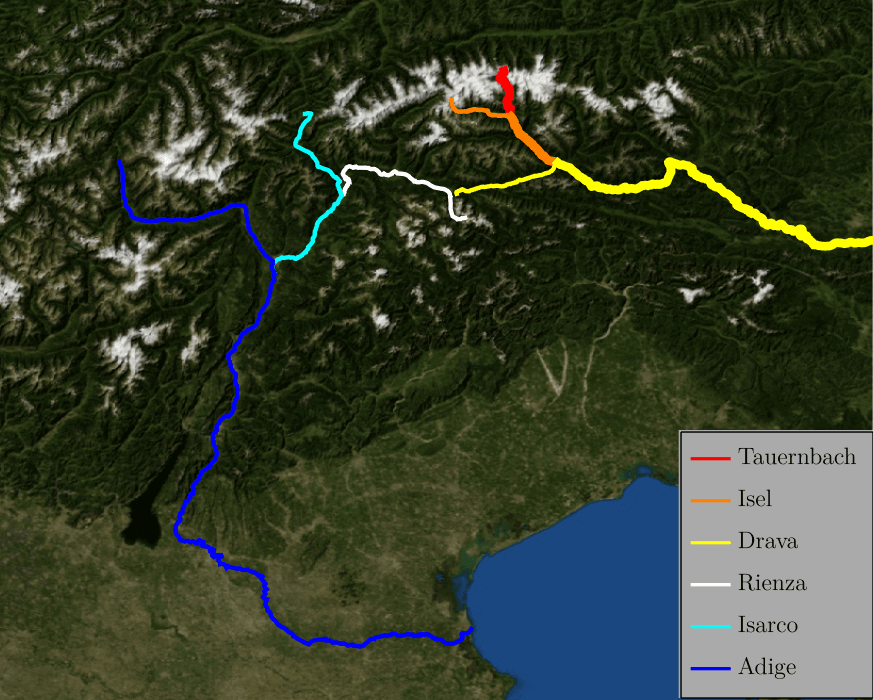}
 \caption{Detail of river courses of  Rienza, Isarco and Adige (white-cyan-blue) and the river Ister by our result, corresponding to the Tauernbach-Isel-Drava-Danube course (red-orange-yellow), plotted thicker.
}
 \label{fig:rieky_global}
 \end{center}
\end{figure}

\vskip 3mm
First of all, at the end of the second sentence, Strabo clearly says that the distance of the Ister source from the recess of Adriatic is about {\it 1000 stadia} (χιλίους σταδίους). The stadium corresponds to 177.7 - 197.3 m, see \cite{Roller} page 33. 
%We also checked other Strabo's geographical descriptions, e.g., distance from Genova to Monaco given in [4-6-1] is 1110 stadia and on the current Google maps it is about 200 km, so 1 stadium corresponds to 0.180 km = 180 m. 
Our mathematical result, presented in the next sections of the paper, shows that the source of the Ister corresponds either to the source of Drava in Sorgenti della Drava in the Val di Pusteria or to the Isel or Tauernbach sources in the High Tauern, see Figure \ref{fig:rieky_global} (this Figure and some further were created by the help of software \cite{mathematica}). We note that the source of the Isel river was in the past considered approximately in the place of the source of nowadays Tauernbach creek and that such an interconnected stream was called the {\it Isola flu} on historical maps \cite{Blaeu, Mercator}, see Figure \ref{fig:isola}. But let us consider the nowadays situation and measure distances from the recess of Adriatic, placed into Monfalcone close to the ancient Roman city Aquileia, to those three river sources. First, let us measure the distance to the Tauernbach source. We placed the Tauernbach source to the highest possible point below the Gro{\ss}venediger (called the Windisch Taurn on Figure \ref{fig:isola}) on its north-east side. The direct geodesic distance (the shortest path on the Earth surface) to the Monfalcone is about 175 km and the geodesic distance through the passable valleys is about 200 km, see Figure \ref{fig:mapa_cesty} upper row. We measured also the distance by using the local roads to the nearest place to the source, in Schildalm, and we got the distance equal to 214 km, see Figure \ref{fig:mapa_cesty} left bottom. With a high probability, such travel may fit very well with the way to that places in the Strabo's time through the Alpine valleys. And we see that all these distances are in perfect agreement with the Strabo's information about approximately 1000 stadia from the recess of Adriatic! When we considered the distances by the local roads to the other two sources, of the Isel in the Hinterbichl at the end of the valley just south of the Gro{\ss}venediger, and of the Drava in the Val di Pusteria we got 216 km and 189 km, see Figure \ref{fig:mapa_cesty}, respectively. They are again very good    estimates of 1000 stadia. This cannot happen in any case when considering the source of Ister in the source of Danube in Schwarzwald, Germany, with a distance to the recess of Adriatic about 640 km, approximately 3500 stadia, highly exceeding the Strabo's {\it Geographica} information. 

Further, in the first part of the second sentence of the cited text, Strabo says that the Ister makes its beginning from the western highest summits (ἄκρων) of the "Germani" people. Indeed, Tauernbach has its spring exactly between the Gro{\ss}venediger and Gro{\ss}glockner, two highest peaks of the High Tauern, see Figure \ref{fig:rieky_global}, and of all the north-eastern Alps, thus it again perfectly fits with the Strabo's description. Here we note that Strabo also explicitly writes that the "germani" means "genuine" (γνησίους) in his time Roman language. This is clearly stated in the translation of the last sentence of section (7-1-2): "γνήσιοι γὰρ οἱ Γερμανοὶ κατὰ τὴν Ῥωμαίων διάλεκτον" by Roller \cite{Roller}, and also from further context of {\it Geographica} it is clear that in general Strabo uses the term "Germani" to denote all genuine people east of the Rhine  (Ῥῆνος) and north and east of the Alps (including the north-eastern Alps themselves). The term "Germani" represents a much wider notion than the 19th century and the nowadays concept of Germans and their language.

All in all, from the above facts it is clear that the location of the Ister source in a close neighbourhood of the High Tauern is the only possibility fulfilling consistently both Strabo's requirements - to be near the highest summits of the Alps east of the Rhine and to be about 1000 stadia from the recess of Adriatic. Moreover, in the first sentence of the cited Greek text above, Strabo writes that the Ister first flows to the south and soon it changes direction from the west to the east up to the Black Sea. This also perfectly corresponds to the Tauernbach-Isel-Drava-Danube course. It flows about 40 km to the south and then in Lienz, after the confluence with the Drava, the course direction is changing from the west to the east. None of these facts can be derived for the Danube river from its source in the Schwarzwald, Germany.

%\begin{figure}[ht]
% \begin{center}
% \includegraphics[width = 100mm]{img/rieky_lokal.png}
% \caption{Detail of river courses of  Rienza and Isarco (white and cyan) and the river Ister by our result, corresponding to the Tauernbach-Isel-Drava-Danube course (red-orange-yellow), plotted thicker.}
% \label{fig:rieky_lokal}
% \end{center}
%\end{figure} 

The second, very important contradiction between placing the Ister source to the Schwarzwald, Germany and Strabo's {\it Geographica} occurs in section (4-6-9). The whole Book 4, Part 6 (4-6) is devoted to a detailed description of the region of Alps from the Savona (Σαβάτα)
in Liguria, Italy up to the Nanos plateau (Ὄκρᾳ)
in Inner Carniola, Slovenia. Strabo's description of the Alps (Ἄλπεις, Ἄλπεια), many times called also Albia (Ἄλβια) in (4-6), follows first the direction from the south to the north, i.e. from Savona up to the Alpine part of the river Rhine.
Then the description turns to the east, see (4-6-8) and (4-6-9), and going very consistently through the countries (even nowadays federal states) and mountainous regions of Ellvettians (Switzerland - Swiss - Helvetica), Boians (Bavaria - Bayern), Rhaetians (Tyrol and South Tyrol), Noricians (Salzburg and Upper Austria), Tauriskians (Styria - Steiermark) up to Karnians (Carinthia - K\"arnten and Carniola - Krain), not far from the recess of Adriatic, above the territory of Karnians, Strabo arrives at the mountainous places where the source of Ister is located. In this part of section (4-6-9), Strabo describes the river Isaras (Ἰσάρας) which after joining with the river Atagis (Ἄταγις)
empties into the Adriatic. Clearly, Atagis corresponds to river Adige and Isaras corresponds to nowadays river Isarco - or most probably - to the course of Isarco continuing in Bressanone upstream by its (larger) tributary Rienza stemming from the Dolomites and flowing through Val di Pusteria, see Figure \ref{fig:rieky_global}. And in these places also Ister (Ἴστρος) takes its beginning, Strabo says explicitly. It is very clear that all the sources, of Tauernbach, Isel or Drava, fulfils this geographic requirement. On the other hand, placing the source of Ister in the Schwarzwald, Germany cannot solve in any way such geographic situation of the Ister flowing from its source to the Black Sea and the neighbouring rivers flowing to the Adriatic. We see that our mathematical result brings the straightforward solution to this tedious and long-lasting controversy which yielded many troubles in translations of {\it Geographica}, even by exchanging the names of the rivers in Section (4-6-9) compared to original, see \cite{MD} where such possible "purposeful" changes were only indicated in brackets in the Latin translation and \cite{Meineke} where such changes of river names were even performed in the Greek text.

Just as a curiosity we mention that our result simply solves also the otherwise unexplainable paradox in the {\it voyage of Argonauts}, where, by Apollonius of Rhodes (3rd century BC) in his epic poem {\it Argonautica}, "the Argonauts had been obliged to abandon their regular course from Colchis homeward, and had gone from Euxine Sea (Black Sea) up the Ister and then passing down the other branch of that river, they had entered into Adriatic" \cite{Grote}.  By our result it is allowed, Argonauts could follow the Danube-Drava upstream up to the Val di Pusteria and continue to the Adriatic by the stream of Rienza-Isarco-Adige. The stream of Rienza is only 3 km away from the source of the Drava in Val di Pusteria where the watershed between the Black Sea and the Adriatic is located, see Figure \ref{fig:rieky_global}. Of course, we do not want to claim that Argonauts made such a voyage :-) but to explain why "such story was accepted even by so able geographer as Eratosthenes who seems to have been a firm believer in the reality of the Argonautic voyage" \cite{Grote}. There is no controversy and it may indicate that Eratosthenes was aware of the watershed in this Alpine region. Moreover, if we adopt a hypothesis that the notion of Ister was evolved in time as Greeks and Romans explored Europe between the Adriatic and the Black Sea from the south, we see also further "earlier" possibility of the Argonauts' voyage, following upstream the Danube-Sava course up to Zelenci and continuing near Tarvisio to the Adriatic by the Fella-Tagliamento stream. 

After the above explanations of solving main geographical controversies appearing when reading {\it Geographica} and its later translations, let us begin our mathematical story. 

It will be explained with all mathematical and computational details in the next Sections 2 and 3. Moreover, since we have to shift the sources of Ister from Schwarzwald to the High Tauern and the upper course of the Ister from the upper Danube to the Drava course, many historical facts from the beginning of the first millennium should be "shifted" as well. That opens many questions which will be discussed in Section 4. The paper will be finished by our conclusions in Section 5.

\section{Locally Affine Globally Laplace (LAGL) Map Transformation}
To register the map $M_1$ to map $M_2$, we use an affine transformation. 
In general, an affine transformation is given by the formula 
\begin{equation}
\mathbf{y}=\bm{A}\,\mathbf{x}+\mathbf{b},
\end{equation}
where $\mathbf{x}=\left(x_1,x_2\right)$ is a point on the map $M_1$, and $\mathbf{y}=\left(y_1,y_2\right)$ is a point on the map $M_2$ and 
\begin{equation}
\bm{A}=\begin{pmatrix}
a_1 & a_2 \\
a_3 & a_4 
\end{pmatrix}
\end{equation} is $2 \times 2$ matrix and 
\begin{equation}
\mathbf{b}= \begin{pmatrix}
b_1 \\
b_2  
\end{pmatrix}
\end{equation}
is a translation vector. For simplicity, we can write the affine transformation as follows
\begin{equation}
\begin{aligned}
y_1=a_1x_1+a_2x_2+b_1\\
y_2=a_3x_1+a_4x_2+b_2.
\end{aligned}
\end{equation}
Our goal is to find the matrix $\bm{A}$ and the vector $\mathbf{b}$ such that 
\begin{equation}
\left|\mathbf{y}-\left(\bm{A}\,\mathbf{x}+\mathbf{b}\right)\right|^2
\label{eq:error_func}
\end{equation}
is minimal for a chosen set of corresponding points $\mathbf{x} \in M_1$ and $\mathbf{y} \in M_2$. Let us have corresponding points $\mathbf{x_1},\dots,\mathbf{x_n}$ and $\mathbf{y_1},\dots,\mathbf{y_n}$, respectively. Such minimization for all corresponding points is equivalent to minimizing
\begin{equation}
\underset{i=1}{\overset{n}\sum}\left(\left(y_{i_1}-a_1x_{i_1}-a_2x_{i_2}-b_1\right)^2+
\left(y_{i_2}-a_3x_{i_1}-a_4x_{i_2}-b_2\right)^2
\right)\label{eq:ucel_fcia}
\end{equation}
with respect to the matrix and translation vector elements. In order to minimize (\ref{eq:ucel_fcia}) we compute the derivatives with respect to $a_1,a_2,b_1,a_3,a_4$ and $b_2$, and set it to $0$. For example, for the element $a_1$ we get
\begin{equation}
\underset{i=1}{\overset{n}\sum}2\left(y_{i_1}-a_1x_{i_1}-a_2x_{i_2}-b_1\right)\left(-x_{i_1}\right)=0 \label{eq:deriv_ucel_fcia}
\end{equation}
and similarly for other elements. Equation (\ref{eq:deriv_ucel_fcia}) can be written in the form
\begin{equation}
a_1\underset{i=1}{\overset{n}\sum}x_{i_1}x_{i_1}+a_2\underset{i=1}{\overset{n}\sum}x_{i_2}x_{i_1}+b_1\underset{i=1}{\overset{n}\sum}x_{i_1}=\underset{i=1}{\overset{n}\sum}y_{i_1}x_{i_1}
\end{equation} 
and from there we can see that for all elements we get the system of linear equations 
\begin{equation}
\scalemath{0.96}{
\begin{bmatrix}
\!\suma x_{i_1}^2\! &\!\!\suma x_{i_1} x_{i_2} &\!\suma x_{i_1} &\!\!\!0 &\!0 &\!0\\
\suma x_{i_1} x_{i_2}\! &\!\!\suma x_{i_2}^2 &\!\suma x_{i_2} &\!\!\!0 &\!0 &\!0\\
\!\suma x_{i_1}\! &\!\!\suma x_{i_2} &\!n &\!\!\!0 &\!0 &\!0 \\
\!0\! &\!\!0 &\!0 &\!\!\!\suma x_{i_1}^2 &\!\suma x_{i_1} x_{i_2} &\!\suma x_{i_1}\\
\!0\! &\!\!0 &\!0 &\!\!\!\suma x_{i_1} x_{i_2} &\!\suma x_{i_2}^2 &\!\suma x_{i_2}\\
\!0\! &\!\!0 &\!0 &\!\!\!\suma x_{i_1} &\!\suma x_{i_2} &\!n
\end{bmatrix}
\!\!\!
\begin{bmatrix}
\vphantom{\suma}a_1\\
\vphantom{\suma}a_2\\
\vphantom{\suma}b_1\\
\vphantom{\suma}a_3\\
\vphantom{\suma}a_4\\
\vphantom{\suma}b_2\\
\end{bmatrix}
\!\!=\!\!
\begin{bmatrix}
\suma y_{i_1} x_{i_1} \\
\suma y_{i_1} x_{i_2} \\
\suma y_{i_1}\\
\suma y_{i_2} x_{i_1} \\
\suma y_{i_2} x_{i_2} \\
\suma y_{i_2}
\end{bmatrix}}
\label{eq:transf_matrix}
\end{equation}
%We can see that above system is split into two independent systems of three equations with three unknowns.

Regarding our application to transform the {\it Strabo's map of the World} to the current map in the WGS, we present now the finding of affine transformations $T^{p_k}$ by (\ref{eq:transf_matrix}) for the selected corresponding points sets $p_k, k=1,2,3$, with the corresponding points given at the Adriatic coast ($k=1$), Greece and Albania region ($k=2$) and at the Black Sea coast ($k=3$), for more geographic information about the corresponding points sets see the beginning of section \ref{sec:application}. The accuracy of the affine transformation $T^{p_k}$ for the corresponding points set $p_j=\{\mathbf{x}_i,\mathbf{y}_i;i=1,...,n_{p_j}\}$ is measured by the mean error
\begin{equation}
\varepsilon_{mean}^{p_j}\left(T^{p_k}\right) = \sqrt{\frac{1}{n_{p_j}}\underset{\mathbf{x}_i \in p_j} \sum D_E\left(\mathbf{y}_i, T^{p_k}\left(\mathbf{x}_i\right)\right)^2},
\end{equation} 
where $D_E\left(\mathbf{y}_i,T^{p_k}\left(\mathbf{x}_i\right)\right)$ is a geodesic distance of point $\mathbf{y}_i$ and transformed point $T^{p_k}\left(\mathbf{x}_i\right)$ computed by the GeographicLib::Geodesic class \cite{geoLib}, We also measure the maximal error of the transformation by
\begin{equation}
\varepsilon_{max}^{p_j}\left(T^{p_k}\right) = \max_{\mathbf{x}_i \in p_j} D_E\left(\mathbf{y}_i,T^{p_k}\left(\mathbf{x}_i\right)\right).
\end{equation}

\begin{table} [!h]
 \caption{The mean errors of the transformation $T^{p_k}$ for the corresponding points sets $p_j$. \label{tab:aff_local_mean}} 
 \begin{center} \footnotesize 
 \begin{tabular}{|c|c|d{3.3}|d{3.3}|d{3.3}|}
  \hline 
  \textit{k} & $p_k$ - region & \multicolumn{1}{c|}{$\varepsilon_{mean}^{p_1\vphantom{^R}}\left(T^{p_k}\right)\;\left[km\right]$} &\multicolumn{1}{c|}{$\varepsilon_{mean}^{p_2\vphantom{^R}}\left(T^{p_k}\right)\;\left[km\right]$} &
  \multicolumn{1}{c|}{$\varepsilon_{mean}^{p_3\vphantom{^R}}\left(T^{p_k}\right)\;\left[km\right]$} \\\hline  \hline
1 & Adriatic coast & 7.527 & 256.856 & 723.619 \\\hline
2 & Greece and Albania & 197.655 & 24.259 & 126.697 \\\hline
3 & Black sea coast & 155.672 & 189.790 & 22.981 \\\hline
\end{tabular} 
 \end{center} 
\end{table} 

\begin{table} [!h]
 \caption{The maximal errors of the transformation $T^{p_k}$ for the corresponding points sets $p_j$. \label{tab:aff_local_max}} 
 \begin{center} \footnotesize 
 \begin{tabular}{|c|c|d{3.3}|d{3.3}|d{3.3}|}
  \hline 
  \textit{k} & $p_k$ - region & \multicolumn{1}{c|}{$\varepsilon_{max}^{p_1\vphantom{^R}}\left(T^{p_k}\right)\;\left[km\right]$} &\multicolumn{1}{c|}{$\varepsilon_{max}^{p_2\vphantom{^R}}\left(T^{p_k}\right)\;\left[km\right]$} &
  \multicolumn{1}{c|}{$\varepsilon_{max}^{p_3\vphantom{^R}}\left(T^{p_k}\right)\;\left[km\right]$}  \\\hline  \hline
1 & Adriatic coast &12.201 & 383.310 & 917.010 \\\hline
2 & Greece and Albania & 235.104 &  32.042 & 178.451\\\hline
3 & Black sea coast & 198.112 & 232.691 & 35.989 \\\hline  
\end{tabular} 
 \end{center} 
\end{table} 

\begin{table} [!h]
 \caption{The mean and maximal errors of the transformation $T^{p_G}$ computed by using the global corresponding points set $p_G = p_1 \cup p_2 \cup p_3$. In the first three lines, the errors of $T^{p_G}$ for the corresponding points just from the sets $p_k=1,2,3$ are presented. Comparing these errors with the errors on diagonals of {\rm Tables \ref{tab:aff_local_mean} - \ref{tab:aff_local_max}} we see the error increase. The last line represents the errors of $T^{p_G}$ for all corresponding points in the global set $p_G$. \label{tab:aff_global} } 
 \begin{center} \footnotesize 
 \begin{tabular}{|c|c|d{3.3}|d{3.3}|}
  \hline 
  \textit{k} & $p_k$ - region &
  \multicolumn{1}{c|}{$\varepsilon_{mean}^{p_k\vphantom{^R}}\left(T^{p_G}\right)\;\left[km\right]$} 
  & \multicolumn{1}{c|}{$\varepsilon_{max}^{p_k}\left(T^{p_G}\right)\;\left[km\right]$} \\\hline  \hline
1 & Adriatic coast & 45.730 & 65.453 \\\hline
2 & Greece and Albania & 62.224 & 92.509 \\\hline
3 & Black Sea coast& 43.519 & 79.303 \\\hline\hline
\textit{G} & global & 52.543 & 92.509 \\\hline
\end{tabular} 
 \end{center} 
\end{table} 

\begin{figure}[!h]
 \begin{center}
 \includegraphics[height = 44.65mm]{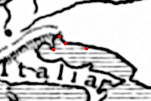}
 \includegraphics[trim={0cm 0 3cm 0},clip,height = 44.65mm]{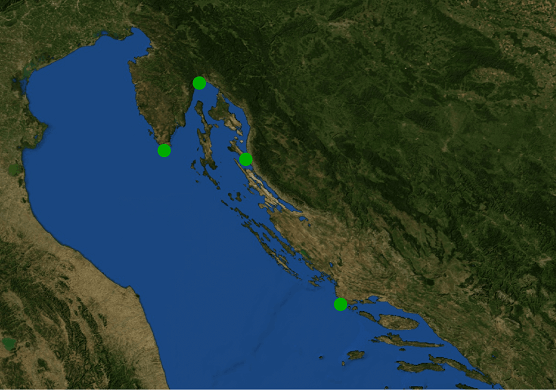}
 \includegraphics[trim={0cm 0cm 0cm 15mm},clip,width=1.\textwidth]{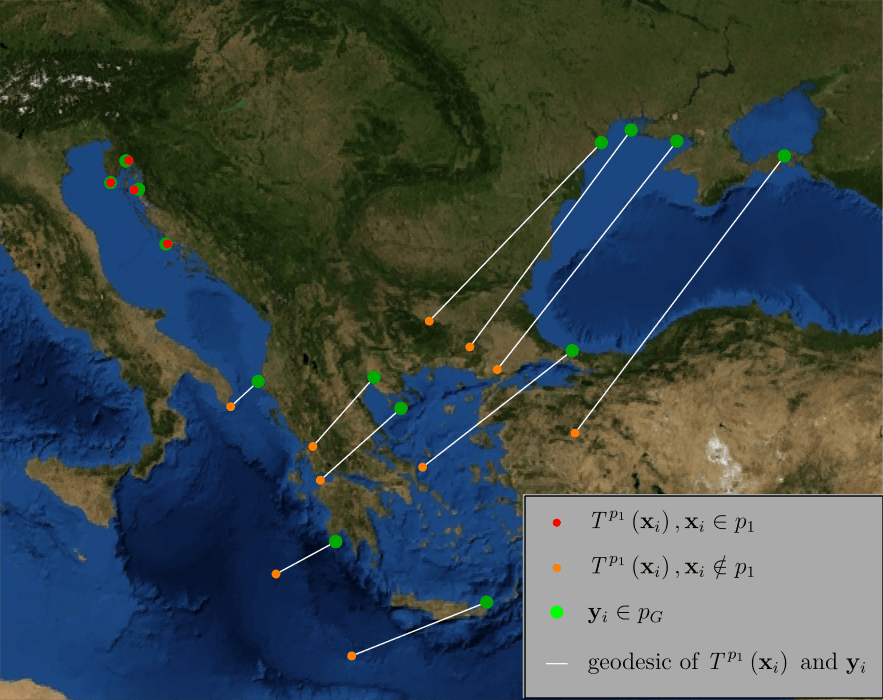}
 \caption{Optimal local affine transformation $T^{p_1}$ obtained by using the corresponding points set $p_1$ given on the Adriatic coast. Top images show the corresponding points, on the left in red on the Strabo's map and on the right in green on the map in WGS. The bottom image shows the transformed points of the set $p_1$ in red, the transformed points of the sets $p_2$ and $p_3$ in orange and their corresponding points in WGS  in green. One can see, that points on Adriatic coast are transformed accurately while the points in other regions are distant from their corresponding points, see also the first rows of {\rm Tables \ref{tab:aff_local_mean} - \ref{tab:aff_local_max}}.
}\label{fig:affineJadran}
 \end{center}
 \end{figure}

\begin{figure}[!h]
 \begin{center}
 \includegraphics[height=44.65mm]{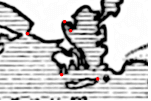}
 \includegraphics[trim={0cm 0 30mm 0},clip,height = 44.65mm]{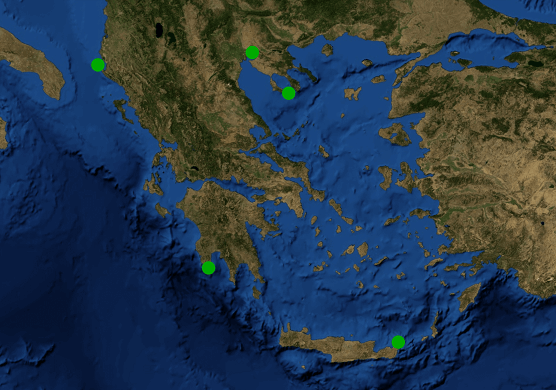}
 \includegraphics[trim={0cm 0cm 0cm 15mm},clip,width = 1.\textwidth]{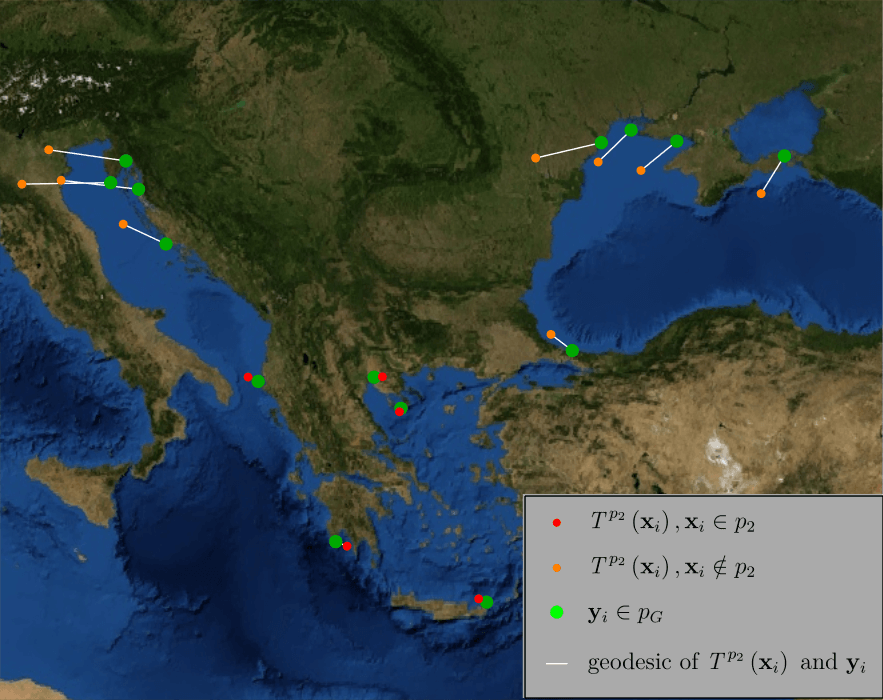}
 \caption{Optimal local affine transformation $T^{p_2}$ obtained by using the corresponding points set $p_2$ given in the Greece and Albania region. Top images show the corresponding points, on the left in red on the Strabo's map and on the right in green on the map in WGS. The bottom image shows the transformed points of the set $p_2$ in red, the transformed points of the sets $p_1$ and $p_3$ in orange and their corresponding points in WGS  in green. One can see, that points in the Greece and Albania region are transformed accurately while the points in other regions are distant from their corresponding points, see also the second rows of {\rm Tables \ref{tab:aff_local_mean} - \ref{tab:aff_local_max}}.}\label{fig:affineGrecko}
 \end{center}
 \end{figure}

\begin{figure}[!h]
 \begin{center}
 \includegraphics[height=44.65mm]{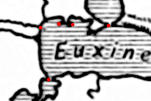}
 \includegraphics[trim={0cm 0 3cm 0},clip,height = 44.65mm]{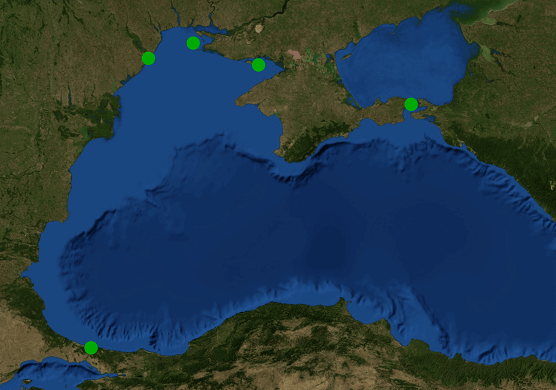}
 \includegraphics[trim={0cm 0cm 0cm 15mm},clip,width = 1.\textwidth]{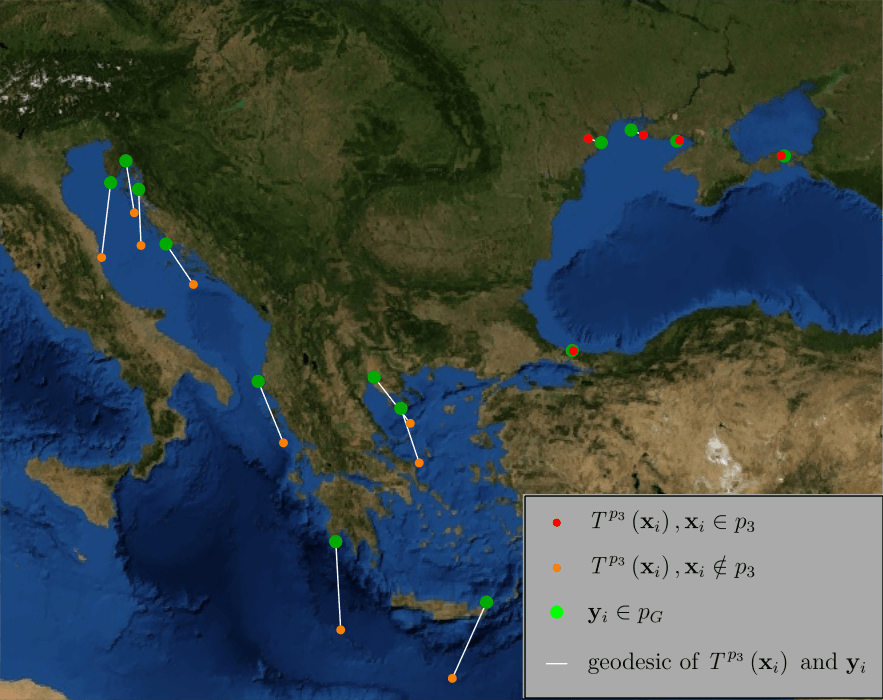}
 \caption{Optimal local affine transformation $T^{p_3}$ obtained by using the corresponding points set $p_3$ given on the Black Sea coast. Top images show the corresponding points, on the left in red on the Strabo's map and on the right in green on the map in WGS. The bottom image shows the transformed points of the set $p_3$ in red, the transformed points of the sets $p_1$ and $p_2$ in orange and their corresponding points in WGS  in green. One can see, that points on Black Sea coast are transformed accurately while the points in other regions are distant from their corresponding points, see also the third rows of {\rm Tables \ref{tab:aff_local_mean} - \ref{tab:aff_local_max}}.
}\label{fig:affineCmore}
 \end{center}
 \end{figure}

\begin{figure}[!h]
 \begin{center}
 \includegraphics[trim={0cm 0cm 0cm 15mm},clip,width = 1.\textwidth]{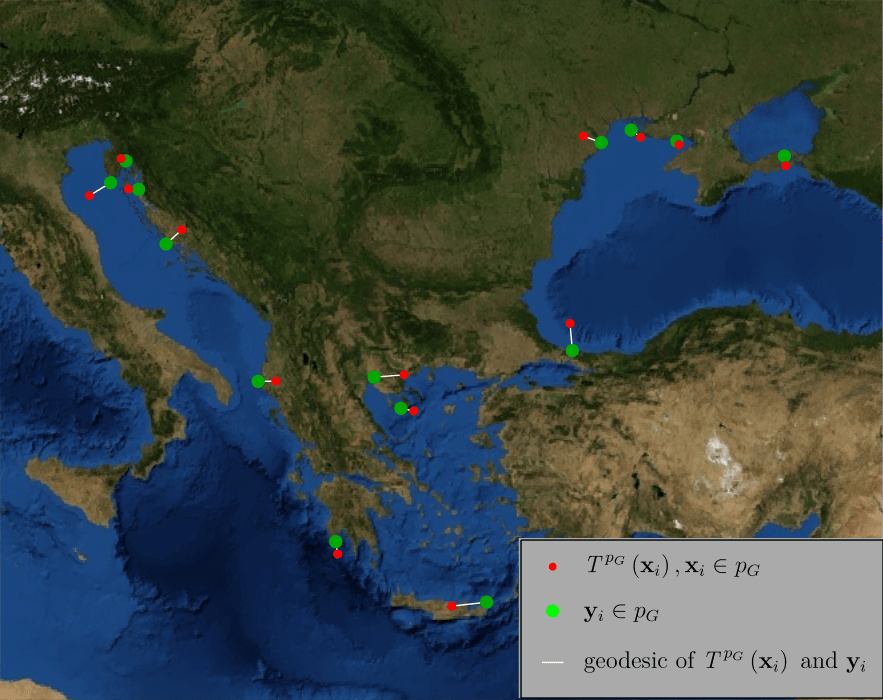}
 \caption{Optimal local affine transformation $T^{p_G}$ obtained by using the corresponding points set $p_G=p_1 \cup p_2 \cup p_3$. The image shows the transformed points from the Strabo's map in red and their corresponding points in WGS in green. One can see that all the points are transformed relatively closely. However, in the previous local transformations, see {\rm Figures \ref{fig:affineJadran} -\ref{fig:affineCmore}}, the points from the sets used for finding the optimal local affine transformations were transformed more accurately, see also {\rm Tables \ref{tab:aff_local_mean} - \ref{tab:aff_global}}.
}\label{fig:affineJGC}
 \end{center}
 \end{figure}

Transformations presented in the Tables \ref{tab:aff_local_mean} - \ref{tab:aff_local_max} comes with acceptable errors for the corresponding points sets $p_k$ if the same corresponding points are used also for finding the transformation $T^{p_k}$, see diagonals of the tables. Unfortunately, the error dramatically rises for the corresponding points not used for finding the optimal transformation, see out of diagonal entries in the tables and Figures \ref{fig:affineJadran} - \ref{fig:affineCmore}. Thus these local transformations are not suitable for transformation of farther points not included in the minimization procedure. Of course, one can create a common global corresponding points set $p_G = p_1 \cup p_2 \cup p_3$ and find by (\ref{eq:transf_matrix}) a common global transformation $T^{p_G}$.  Using such optimal affine transformation $T^{p_G}$ we get visually better results, see Fig \ref{fig:affineJGC}. However, the errors $\varepsilon_{mean}^{p_k}\left(T^{p_G}\right)$ and $\varepsilon_{max}^{p_k}\left(T^{p_G}\right)$, $k=1,2,3$, see Table \ref{tab:aff_global}, have significantly increased compared to the corresponding errors on the diagonals of Tables \ref{tab:aff_local_mean} - \ref{tab:aff_local_max}. Also visually we see quite large differences compared to the local transformations. For example, for the Adriatic coast region the mean error rises from 7.5 km to 45.7 km and maximal error from 12.2 km to 65.4 km which is expressed visually in Figures \ref{fig:affineJadran} and \ref{fig:affineJGC}. A similar result is observed near Istanbul and in other localities as well. Since we want to keep the accuracy of the local affine transformations in the neighbourhood of selected polygonal regions and not to pollute the transformations by the globally increasing errors we come with the following idea.
 
\clearpage

\clearpage
 
By using the points $\mathbf{x}_i \in p_k$ we create a polygon denoted again without any confusion by $p_k$, and it is done for all $k=1,...,n_p$ where $n_p$ is the number of corresponding points sets and the created polygons as well. In all points of the map $M_1$, which are inside of every polygon $p_k$, we set the parameters of affine transformation to the values of locally optimal affine transformation $T^{p_k}$ found by (\ref{eq:transf_matrix}) using the corresponding points set $p_k$. For all other points of the map $M_1$ we use an interpolation/extrapolation approach. In 1D case, when we want to interpolate between two given function values, the natural approach is to use the linear interpolation, i.e. to connect two given values by a straight line. However, it is not so straightforward in higher dimensions. Fortunately, there exists an analogy of linear interpolation in higher dimensions given by a solution of the Laplace equation with Dirichlet boundary conditions given on the boundary of the interpolation domain. It is widely used in data processing such as filling the missing parts of photographs or other image inpainting problems. However, in our case, we have not only to interpolate the values to the regions between the polygons but also to extrapolate these values to the whole map $M_1$. To that goal, we suggest the following mathematical model. Let us consider the Laplace equation
\begin{equation}
- \Delta u\left(\mathbf{x}\right) =0
\label{eq:laplace}
\end{equation}
where solution $u$ represents the transformation matrix and translation vector elements $a_1,a_2,b_1,a_3,a_4, b_2$, together with the Dirichlet conditions prescribed in the polygons $p_k$. The Laplace equation with such Dirichlet conditions is solved in a domain $\Omega$ and the zero Neumann boundary conditions are prescribed on its boundary $\partial\Omega$. The domain $\Omega$ is chosen as a rectangular subset of the map $M_1$, see e.g. the pictures in Figure \ref{fig:laplace_interpolation2} where we choose as domain $\Omega$ a rectangle surrounding Europe on the Strabo's map of the World. We note that the minus sign in the equation (\ref{eq:laplace}) is chosen to have operator on the left hand side positive and arising matrix of the system then positive definite.

\begin{figure}
 \begin{center}
 \includegraphics[width = 1.\textwidth]{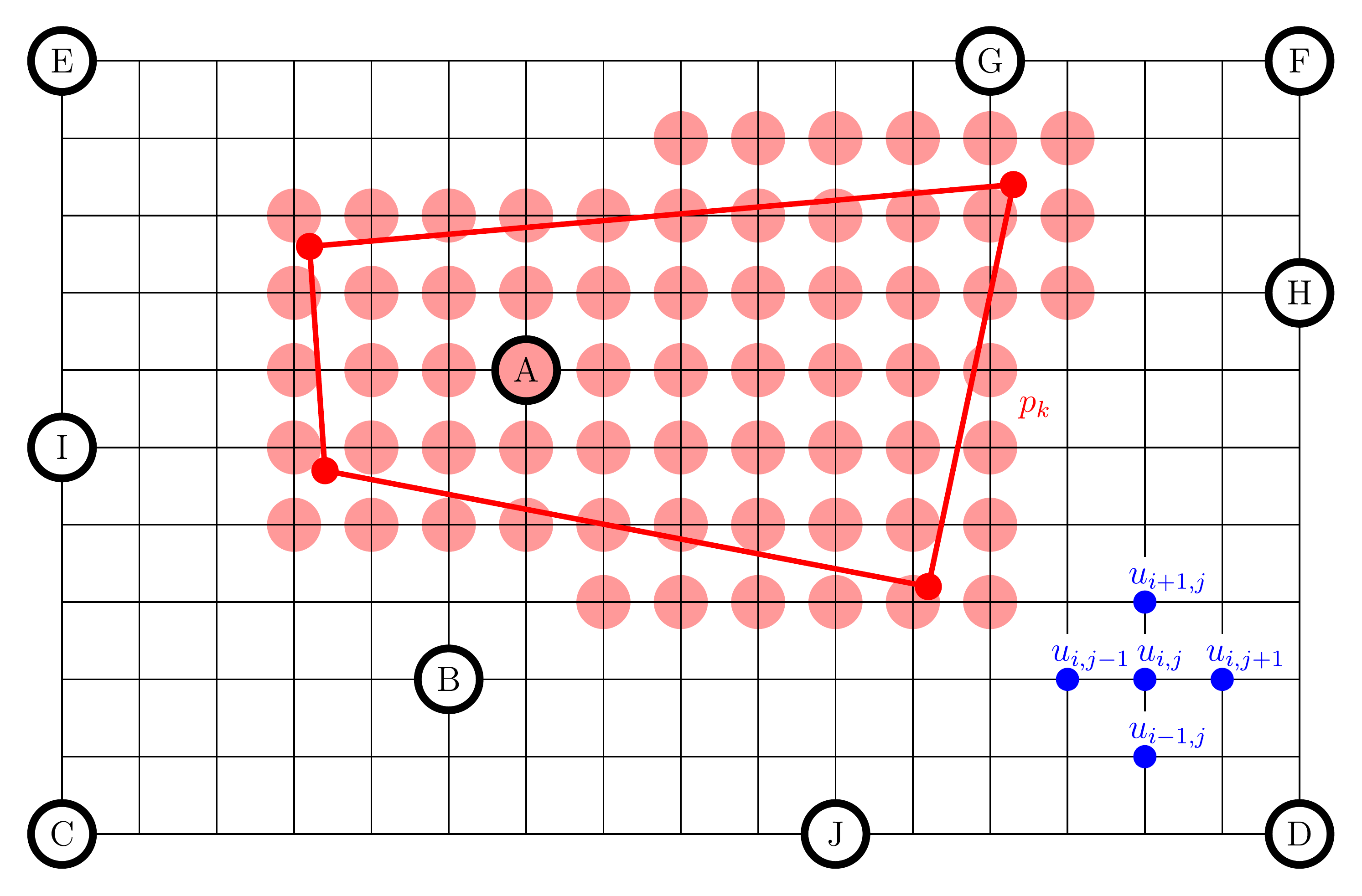}
 \caption{Illustration of the discretization of the computational domain $\Omega$ and approximation of our mathematical model in cases $\rm{A}-\rm{J}$. By pink color we plot the grid nodes in the outer discrete envelope $E(p_k)$ of the polygon $p_k$ (plotted in red) where we consider Dirichlet conditions. In case {\rm B} we consider the standard approximation of the Laplace equation while in cases $\rm{C}-\rm{J}$ we consider its adjustment at the boundary $\partial \Omega$. }
 \label{fig:raster}
 \end{center}
 \end{figure}

To discretize the partial differential equation (\ref{eq:laplace}) in the domain $\Omega$ we use the finite difference method on a uniform grid with the grid size 1, see Fig. \ref{fig:raster}. The grid nodes $\mathbf{x}_{i,j}, i=1,\dots,N_1, j=1,\dots,N_2$ correspond to centers of pixels in the map $M_1$. The Dirichlet conditions are prescribed in $E(p_k)$, the outer discrete envelope of the polygons $p_k$, see Fig. \ref{fig:raster} case A.
In other grid nodes except the boundary, case B in Fig. \ref{fig:raster}, the Laplace operator in equation (\ref{eq:laplace}) is approximated by
\begin{align}
\Delta u\left(\mathbf{x}\right) & = \frac{\partial^2u}{\partial{x_1}^2 }\left(\mathbf{x}\right)+ \frac{\partial^2u}{\partial{x_2}^2 }\left(\mathbf{x}\right)\\
& \approx  u_{i-1,j}-2u_{i,j}+u_{i+1,j}
	+  u_{i,j-1}-2u_{i,j}+u_{i,j+1}
	\\
& =  u_{i-1,j}+u_{i+1,j}
	-4u_{i,j}
	+  u_{i,j-1}+u_{i,j+1},
\end{align}
where $ u_{i,j}$ is an approximate value of $u$ at the grid node $\mathbf{x}_{i,j}$. Such approximation needs to be adjusted for the grid points on the boundary $\partial\Omega$. To approximate the zero Neumann boundary condition we use the reflection of values along the boundary, e.g. in case J in Fig. \ref{fig:raster} we set $u_{i-1,j}=u_{i+1,j}$ which leads to the following approximation on that part of the boundary
\begin{align}
\Delta u\left(\mathbf{x}\right) &\approx  u_{i+1,j}-2u_{i,j}+u_{i+1,j}
	+  u_{i,j-1}-2u_{i,j}+u_{i,j+1}
	\\
& =  2u_{i+1,j}
	-4u_{i,j}
	+  u_{i,j-1}+u_{i,j+1}
\end{align}
and it is done similarly for other boundary grid nodes. 

We summarize all discrete equations, for cases A-J, representing the numerical discretization of our model for the LAGL map transformation as follows:
\vskip 2mm
\begin{itemize}
\item[A:] $\mathbf{x}_{i,j} \in E(p_k) : u_{i,j} = v\left(p_k\right) $, where $v$ is any of the element of $\bm{A}$ and $\mathbf{b}$ given by locally optimal affine transformation $T^{p_k}$ found by (\ref{eq:transf_matrix}), 
\item[B:] $\mathbf{x}_{i,j} \not\in E(p_k) \;\land\; \mathbf{x}_{i,j} \not\in \partial\Omega  : -u_{i-1,j}-u_{i,j-1} + 4u_{i,j} - u_{i+1,j} - u_{i,j+1} = 0$
\item[C:] $\mathbf{x}_{i,j} \not\in E(p_k) \;\land\;
 i=1, j=1
 :   4u_{i,j}-2u_{i+1,j} - 2u_{i,j+1}=0$
\item[D:] $\mathbf{x}_{i,j} \not\in E(p_k) \;\land\;
 i=1, j=N_2
 : - 2u_{i,j-1}+ 4u_{i,j}  -2u_{i+1,j}=0$
\item[E:] $\mathbf{x}_{i,j} \not\in E(p_k) \;\land\;
 i=N_1, j=1
 :  -2u_{i-1,j}+ 4u_{i,j} - 2u_{i,j+1}=0$
\item[F:] $\mathbf{x}_{i,j} \not\in E(p_k) \;\land\;
 i=N_1, j=N_2
 :  -2u_{i-1,j} - 2u_{i,j-1}+ 4u_{i,j}=0$
\item[G:]  $\mathbf{x}_{i,j} \not\in E(p_k) \;\land\;
 i=N_1, j=2,\dots,N_2-1
 :  -2u_{i-1,j}-u_{i,j-1} + 4u_{i,j} - u_{i,j+1}=0$
\item[H:]  $\mathbf{x}_{i,j} \not\in E(p_k) \;\land\;
i=2,\dots,N_1-1, j=N_2
 :  -u_{i-1,j} - 2u_{i,j-1}+ 4u_{i,j}-u_{i+1,j}=0$
\item[I:]  $\mathbf{x}_{i,j} \not\in E(p_k) \;\land\;
 i=2,\dots,N_1-1, j=1
 :  -u_{i-1,j} + 4u_{i,j}-u_{i+1,j} - 2u_{i,j+1}=0$
\item[J:]  $\mathbf{x}_{i,j} \not\in E(p_k) \;\land\;
 i=1, j=2,\dots,N_2-1
 :  -u_{i,j-1} + 4u_{i,j} -2u_{i+1,j}- u_{i,j+1}=0$
\end{itemize}
\vskip 2mm

\begin{figure}
 \begin{center}
 \includegraphics[trim={0cm 0cm 0cm 15mm},clip,width = 1.\textwidth]{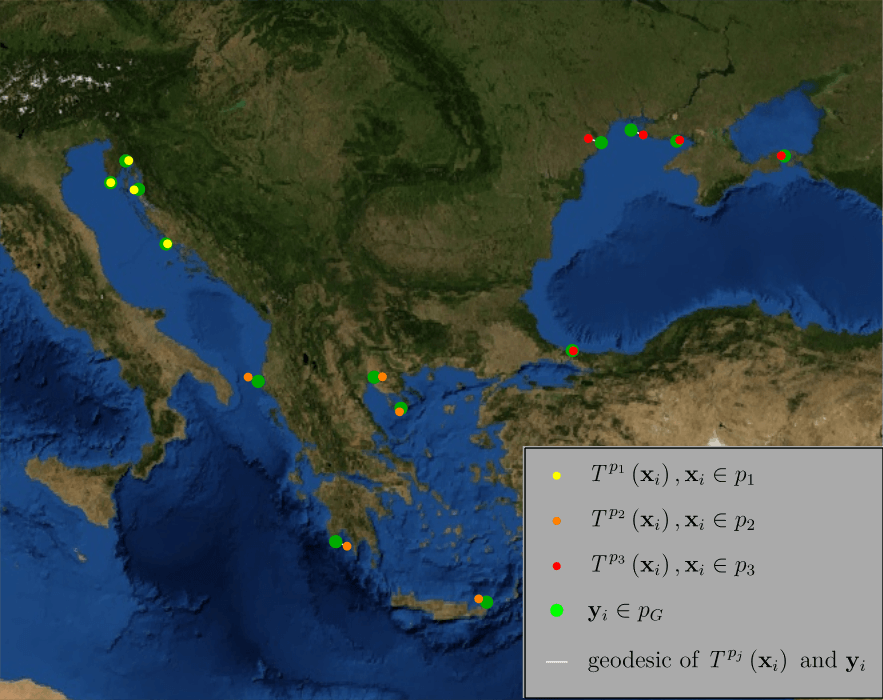}
 \caption{Example of transformation found by using the LAGL method. Top image illustrates the result for the corresponding points sets $p_1$, $p_2$ and $p_3$ which are transformed as accurately as by their locally optimal affine transformations $T^{p_1}$, $T^{p_2}$ and $T^{p_3}$, see {\rm Figures \ref{fig:affineJadran} - \ref{fig:affineCmore}}. 
}
\label{fig:laplace_interpolation1}
 \end{center}
 \end{figure}

\begin{figure}
 \begin{center}
 \includegraphics[width = 1.\textwidth]{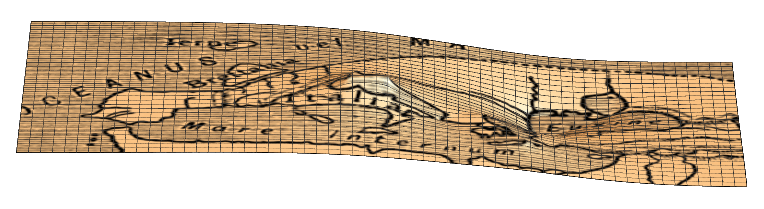}
 \includegraphics[width = 1.\textwidth]{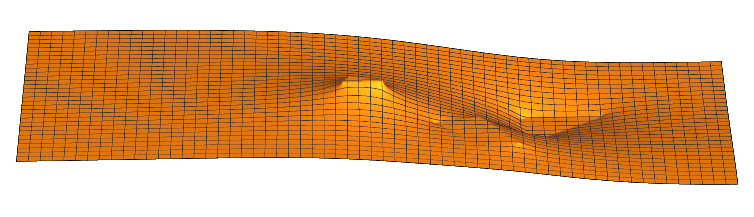}
 \caption{Example of transformation found by using the LAGL method.  
The images show one of the transformation parameters computed by the LAGL method and presented in {\rm Figure \ref{fig:laplace_interpolation1}}. The upper image is plotted with the texture of the Strabo's map of the World, the lower image is plotted without the texture emphasizing the smooth transition of locally optimal affine transformations.}
\label{fig:laplace_interpolation2}
 \end{center}
 \end{figure}

For solving the above system of equations we use Eigen::SparseLU class of \cite{eigen} which solves the linear systems with sparse matrices directly and efficiently by the LU decomposition. 
Solution of the linear system of equations gives us the desired result, smoothly interpolated/extrapolated locally optimal affine transformations from polygons $p_k$ to every point of the (historical) map $M_1$. As we see in Figure \ref{fig:laplace_interpolation1} we reproduce the optimal errors (diagonals in Tables  \ref{tab:aff_local_mean} - \ref{tab:aff_local_max}) of all local affine transformations thanks to the Dirichlet conditions in polygons $p_k$. These optimal values are smoothly interpolated/extrapolated to the whole historical map $M_1$, as seen in Fig. \ref{fig:laplace_interpolation2}. Such locally optimal smoothly varying transformation can be used to transform any point from the historical map $M_1$ to the current map $M_2$.

%\clearpage

\section{Strabo's Ister (Ἴστρος) transformation\label{sec:application}} Now we present the application of the LAGL method to the transformation of the river Ister from the {\it Strabo's map of the World} to the current map in WGS. The corresponding points on Strabo's map are chosen on the east, south and west directions from the Ister river on the {\it Strabo's map of the World}  and in such a way that they are reliably identifiable also on the current map of the world. It is clearly observable that the {\it Strabo's map of the World} by Karl M\"uller is less realistic in the northwestern and northern directions from the Alps than in other directions. But this fact correctly respects the uncertainty of Strabo’s description of this part of Europe and Karl M\"uller's map does not introduce any artificial information which can be found in others 19th-century map reconstructions of the {\it Geographica}. While southern Europe, west of the Alps up to Lyon (Λούγδουνον) and east of the Alps up to the Ister, is described in {\it Geographica} with a lot of quantitative geographic details including distances, the northern part contains almost no such information. The almost only quantitative information is what Strabo says about the length of the Rhine: "Asinius (historian) says that length of the river (Rhine) is six thousand stadia, but it is not, for it would be only a little more than half that in a straight line, and adding a thousand would be sufficient for bending.", see \cite{Roller}, section (4-3-3). The estimate of six thousand stadia (by Asinius) is however a right one since the current length of the Rhine is about 1230 km, but Strabo does not accept that and almost halved it. This shows how uncertain was the knowledge about the river Rhine on the northern planes in the time of Strabo and due to that also on the Karl M\"uller’s map it is much shorter than in the reality. About the northern part of Europe Strabo also writes that the region "near the Ocean is totally unknown to us", see \cite{Roller}, section (7-2-4), and by the "Ocean", he means the nowadays Northern and Baltic Seas. That is the reason why Karl Muller plotted the northern border of Europe (and of the World) by the dashed line. Due to the above reasons, we are not allowed to choose the corresponding points north-west and north of the Ister. 

Having all these requirements in mind, we defined the first three sets of corresponding points on the Mediterranean and the Black Sea costs - in the south-west, south and east directions from the Strabo's Ister, see Figures \ref{fig:affineJadran} - \ref{fig:affineJGC} and Figure \ref{fig:LAGL_JGC}. Since the Mediterranean world up to the Ister river and along it as well as the regions along the Black Sea were already well known for Greeks and Romans in the time of Strabo, the chosen points are reliable regarding location, distances and directions. The first three corresponding points sets defined below fulfils all the above assumptions, but we created also the fourth one containing the points on the Italian coastline and in the Alps. The points in the Alps - the sources of Rhone (Ῥοδανὸς) and Rhine rivers seen on the Strabo's map - have some uncertainty to be identified correctly on the current map. We put the "source" of the Rhine below the Rheinwaldhorn, since in section (4-6-6) Strabo writes that it flows from the mount Adula (Ἀδούλας), and the "source" of Rhona to the Lake Geneva (lac L\'eman) because in (4-6-6) Strabo writes that it bears from Λημέννα λίμνη. We use this fourth corresponding points set mainly to show that considering the western Alpine region, although a bit uncertain, does not change the results significantly. At this place, we would also like to note (probably the only one) discrepancy of the Karl M\"uller’s map and Strabo's descriptions in {\it Geographica}. It is the sketch of the Alps, namely, the "second ridge" should encompass the sources of all the Ister, Rhein and the Rhone because all these river sources are within the Alps by Strabo. We think that Karl M\"uller’s sketch of the "second ridge" is related to the 19th – 20th-century standard assumption of the correspondence of the Strabo's Ister and nowadays Danube rivers.
 
Here are the corresponding points sets:

\vskip 2mm
Adriatic coast region:
\begin{itemize}
\item south of Istria (Premantura)
\item Opatija
\item Jablanac
\item Split area (Ra\v{z}anj)
\end{itemize}

Black Sea coast region:
\begin{itemize}
\item Istanbul (north of Bosporus)
\item mouth of Dniester (Zatoka)
\item cape of Tendrivska gulf
\item cape of Dzharylhatska gulf
\item Ker\v{c}
\end{itemize}

Greece and Albania region:
\begin{itemize}
\item Vlor{\"e}
\item Koufasaratsia
\item north-east cape of Crete (Kyriamadi)
\item cape of Kassandra peninsula
\item Thessaloniki
\end{itemize}

Alps and Italy region:
\begin{itemize}
\item Ancona (Conero)
\item cape south of Venice
\item Trieste
\item "source" of Rhine (Rheinwaldhorn) 
\item "source" of Rhone (Lake Geneva)

\end{itemize}
\vskip 2mm

In numerical experiments presented in this section we vary following combinations of the above-defined regions: 
\begin{itemize}
\item[] Experiment 3.1: Adriatic coast, Black Sea coast regions,
\item[] Experiment 3.2: Adriatic coast, Greece and Albania, Black Sea coast regions,
\item[] Experiment 3.3: Adriatic coast, Black Sea coast, Alps and Italy regions.
\end{itemize}
\vskip 2mm

We will use the following abbreviations:
\begin{itemize}
\item D - Danube
\item D1 - the Danube from the source up to the confluence with the Drava
\item D2 - the Danube from the confluence with the Drava up to outlet to the Black Sea
\item DD - courses of Drava and Danube interconnected
\item TID - courses of Tauernbach, Isel and Drava up to the confluence with the Danube
\item TIDD - courses of Tauernbach, Isel, Drava and Danube interconnected
\item I1 - transformed Ister from the source up to the intersection with the Danube
\item I2 - transformed Ister from the intersection with Danube up to outlets to the Black Sea
\end{itemize}

\begin{figure}[!hb]
 \begin{center}
 \includegraphics[width = 1.\textwidth]{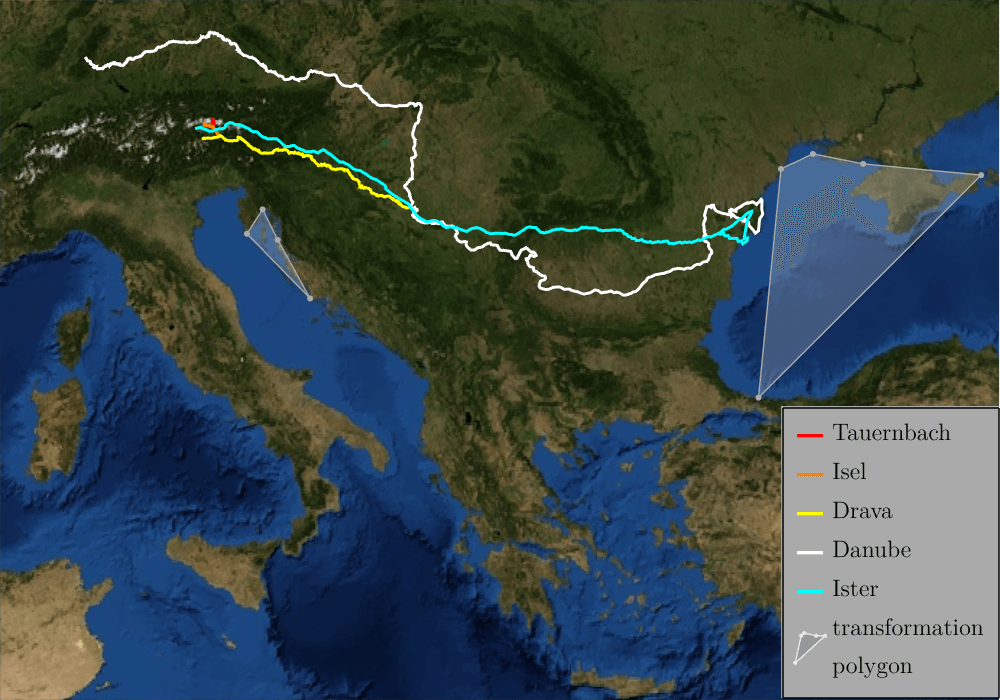}
 \caption{Visual comparison of current rivers (various colors) and river Ister (cyan) transformed using Adriatic and Black Sea coast regions, the polygons used for the transformations are highlighted in grey. \label{fig:LAGL_JC}}
 \end{center}
 \end{figure}

\begin{figure}[!hb]
 \begin{center}
 \includegraphics[width = 1.\textwidth]{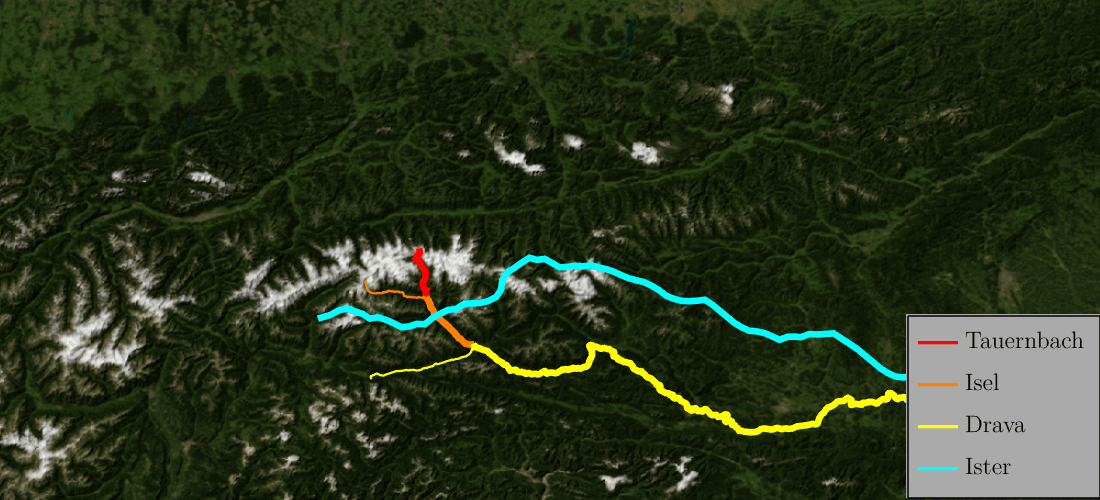}
 \caption{Detail of real rivers source (various colors) compared to the source of transformed river Ister (cyan). We plot the course Tauernabach-Isel-Drava thicker since we consider it to be Ister.
}
\label{fig:LAGL_JC_local}
 \end{center}
 \end{figure}

\begin{figure}[!hb]
 \begin{center}
 \includegraphics[width = 1.\textwidth]{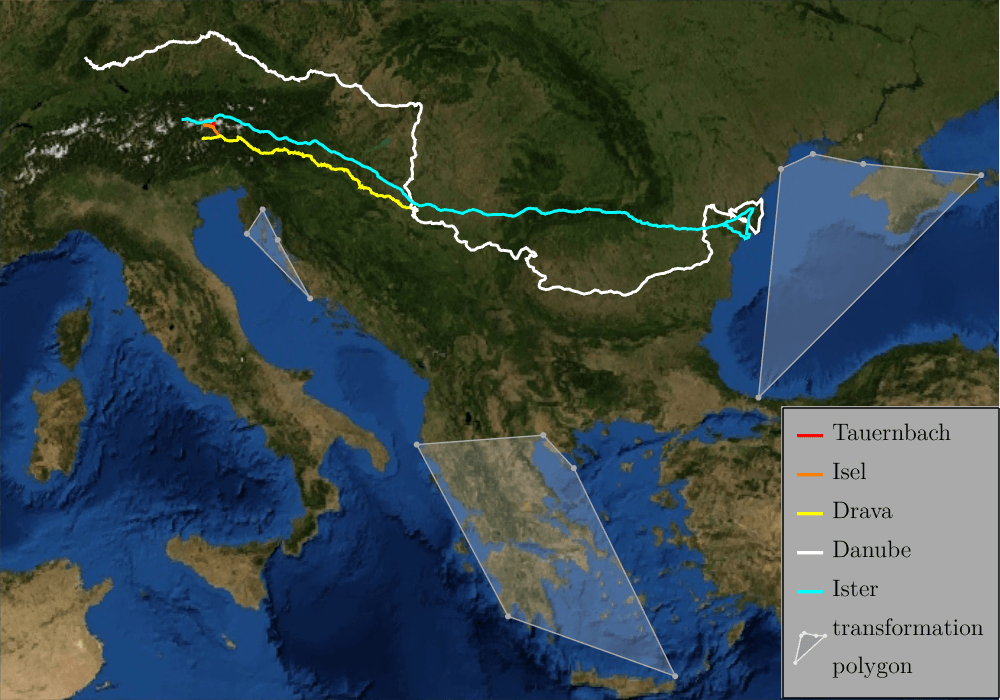}
 \caption{Visual comparison of current rivers (various colors) and river Ister (cyan) transformed using Adriatic, Black sea and Greece and Albania regions, the polygons used for the transformations are highlighted in grey.}
 \label{fig:LAGL_JGC}
 \end{center}
 \end{figure}

\begin{figure}[!hb]
 \begin{center}
 \includegraphics[width = 1.\textwidth]{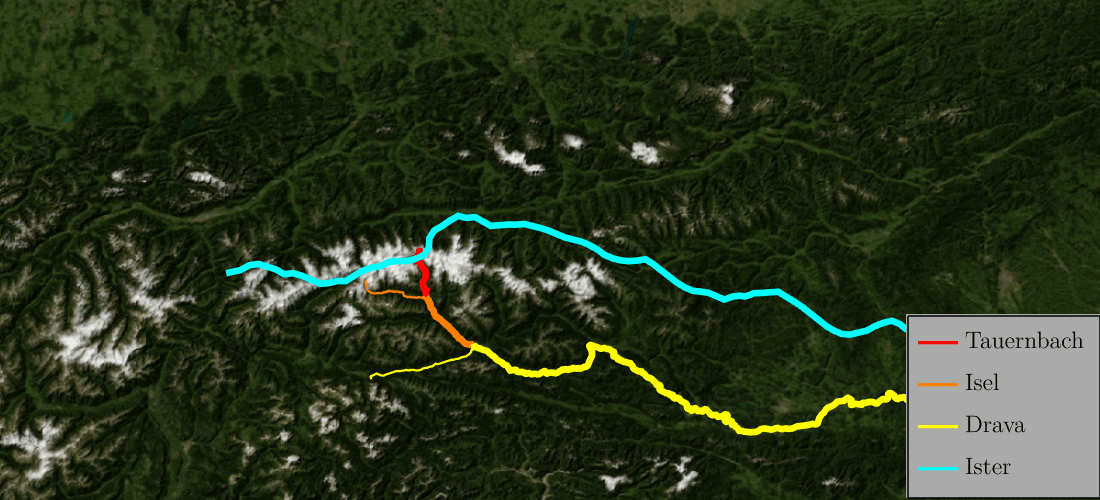}
 \caption{Detail of real rivers source (various colors) compared to the source of transformed river Ister (cyan). We plot the course Tauernabach-Isel-Drava thicker since we consider it to be Ister.}
 \label{fig:LAGL_JGC_local}
 \end{center}
 \end{figure}

\begin{figure}[h!]
 \begin{center}
 \includegraphics[width = 1.\textwidth]{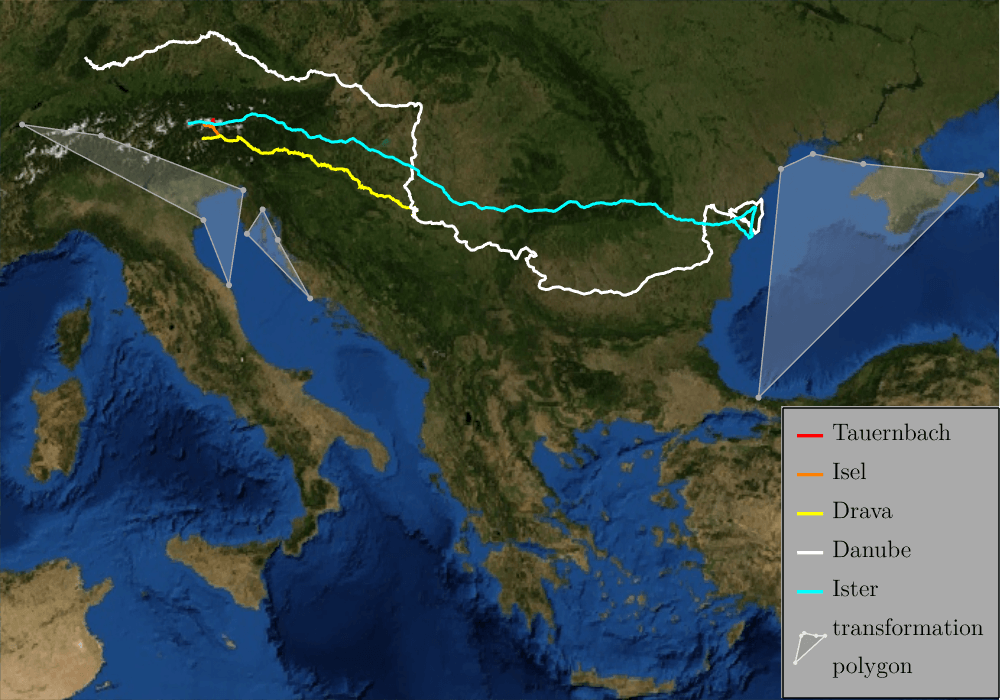}
 \caption{Visual comparison of current rivers (various colors) and river Ister (cyan) transformed using Adriatic, Black sea and Alps and Italy regions, the polygons used for the transformations are highlighted in grey.}
 \label{fig:LAGL_AJC}
 \end{center}
 \end{figure}

\begin{figure}[h!]
 \begin{center}
 \includegraphics[width = 1.\textwidth]{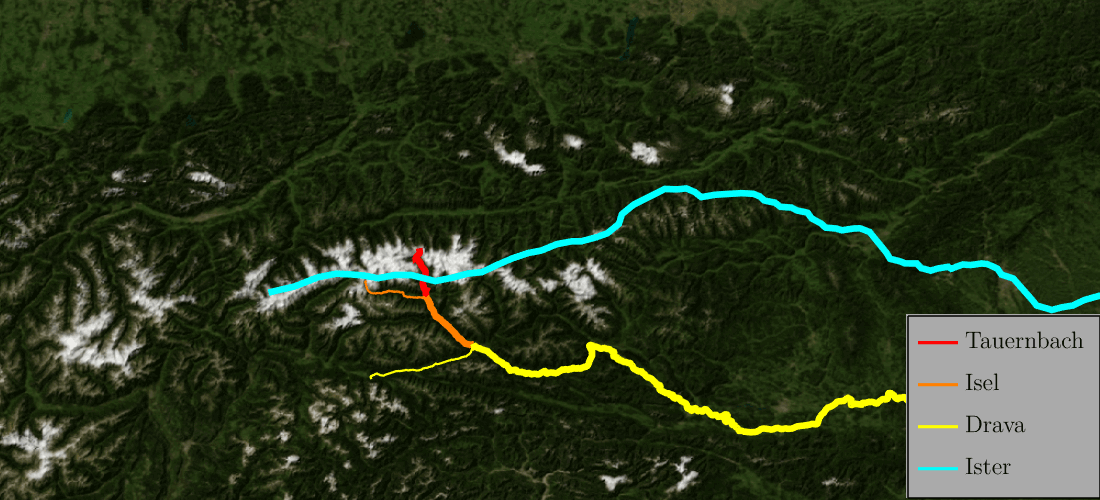}
 \caption{Detail of real rivers source (various colors) compared to the source of transformed river Ister (cyan). We plot the course Tauernabach-Isel-Drava thicker since we consider it to be Ister.
 \label{fig:LAGL_AJC_detail}}
 \end{center}
 \end{figure}

\vskip 2mm
To evaluate the results of transformations we compute the maximal and the mean Hausdorff distances (defined below) of the two discrete curves - one representing the real river course, precisely digitized on the current map, and one representing the discrete transformed Ister from the Strabo's map to the current map. To get the transformed Ister, first, the Ister on the Strabo's map, see Figure \ref{fig:mapa_strabo}, was digitized to contiguous pixel set and then every center of the pixel was transformed to the current map by the LAGL map transformation. All necessary distances on the current map in WGS are computed by means of the GeographicLib::Geodesic class \cite{geoLib}. In Figures \ref{fig:LAGL_JC} - \ref{fig:LAGL_AJC_detail}, the cyan curve represents always the transformed Strabo's Ister while the white curve represents the Danube river, yellow curve the Drava river, the orange curve represents the Isel river and red curve the Tauernbach. We also measure the length of two curves matching in a prescribed narrow band by the so-called matching length defined below. 

Let us have a discrete curve $A=\left\{\mathbf{a}_1,\dots,\mathbf{a}_{n_A}\right\}$. By using the points $\mathbf{a}_i, i=1,\dots, n_A$ we create the piecewise linear segments
$\widehat{A}=\left\{\widehat{\mathbf{a}}_1,\dots,\widehat{\mathbf{a}}_{n_A}\right\}$ 
as follows
\begin{align}
&\widehat{\mathbf{a}}_1=\overline{\mathbf{a}_1,\frac{\mathbf{a}_1+\mathbf{a}_2}{2}}, \\
&\widehat{\mathbf{a}}_i=\overline{\frac{\mathbf{a}_{i-1}+\mathbf{a}_{i}}{2},\mathbf{a}_i}
\cup
\overline{\mathbf{a}_i,\frac{\mathbf{a}_{i}+\mathbf{a}_{i+1}}{2}}, \;\;\;i=2,\dots,n_A-1\\
&\widehat{\mathbf{a}}_{n_A}=\overline{\frac{\mathbf{a}_{n_A-1}+\mathbf{a}_{n_A}}{2},\mathbf{a}_{n_A}},
\end{align}
where $\overline{\mathbf{u},\mathbf{v}}$ represents the line segment connecting points $\mathbf{u}$ and $\mathbf{v}$. Let $\widehat{a}_i=\left|\widehat{\mathbf{a}}_i\right|$ be the piecewise linear segment length and let $L_A$ be the length of the overall discrete curve $A$ given by the sum of the segments length.

The so-called directed mean Hausdorff distance $\overline{\delta}_{H}(A,B) $ of two discrete curves $A=\left\{\mathbf{a}_1,\dots,\mathbf{a}_{n_A}\right\}$ and $B=\left\{\mathbf{b}_1,\dots,\mathbf{b}_{n_B}\right\}$ with segments
$\widehat{A}=\left\{\widehat{\mathbf{a}}_1,\dots,\widehat{\mathbf{a}}_{n_A}\right\}$ and $\widehat{B}=\left\{\widehat{\mathbf{b}}_1,\dots,\widehat{\mathbf{b}}_{n_B}\right\}$ is given by 
\begin{equation}
\overline{\delta}_{H}(A,B)  =\frac{1}{L_A}
\underset{i=1}{\overset{n_A}{\sum}}\widehat{a}_i\underset{\widehat{\mathbf{b}}_j\in \widehat{B}}\min \;D_E(\mathbf{a}_i,\widehat{\mathbf{b}}_j),
\end{equation}
where $D_E(\mathbf{a}_i,\widehat{\mathbf{b}}_j)$ is the geodesic distance of the point $\mathbf{a}_i$ and the segment $\widehat{\mathbf{b}}_j\in \widehat{B}$. 
Then the {\bf mean Hausdorff distance} $\overline{\text{d}}_{H}(A,B) $ is given by the following formula
\begin{equation}
\overline{\text{d}}_{H}(A,B)  = \frac{L_A\overline{\delta}_{H}(A,B)+L_B\overline{\delta}_{H}(B,A)}{L_A+L_B}.
\end{equation}
The so-called directed maximal Hausdorff distance $\delta_{H}(A,B)$ is given by
\begin{equation}
\delta_{H}(A,B)  =  
\sup_{\mathbf{a}_i\in A\vphantom{\widehat{\mathbf{b}}}} 
\Inf
_{\widehat{\mathbf{b}}_j\in \widehat{B}} 
   D_E\left(\mathbf{a}_i,\widehat{\mathbf{b}}_j\right)
\end{equation}
and then the {\bf maximal Hausdorff distance} $\text{d}_{H}(A,B)$ is given by
\begin{equation}
\text{d}_{H}(A,B)  = 
\max
\left\{ 
\delta_{H}\left(A,B\right),\delta_{H}\left(B,A\right)\right\}.
\end{equation}
Finally, the sum of length of segments $\widehat{a}_i \in \widehat{A}$  within the given threshold distance $d_t$ from the segments of the curve $B$ gives the {\bf matching length} $L_{m}\left(A,B,d_t\right)$ by
\begin{equation}
L_{m}\left(A,B,d_t\right) = 
\underset{\widehat{\mathbf{a}}_i\in \widehat{A}}{\sum}\widehat{a}_i, \quad\quad \min_{\widehat{\mathbf{b}}_j \in \widehat{B}} D_E(\mathbf{a}_i,\widehat{\mathbf{b}}_j)<d_t.\\
\end{equation}

In the following Figures and Tables, we present the results of Experiments 3.1 - 3.3 representing three different LAGL transformations of the Strabo's Ister to the current map in WGS. In Figures, we evaluate the results visually and in Tables quantitatively. 

In Figures 3.x.1 (x=1,2,3) we visualize the polygons used for finding the locally optimal affine transformations used in LAGL method, the transformed river Ister (cyan) and the rivers Danube, Drava, Isel and Tauernbach (various colors). In all these Figures the Ister in its upper course is really close to the Drava/Tauernbach-Isel-Drava(TID) river courses. 

In Figures 3.x.2 we compare visually the river sources with the source of the transformed Ister. We see that all sources of Drava, Isel and Tauernbach are geographically very close to the source of the transformed Ister. 

Tables 3.x.1 show that the sources of Danube and transformed Ister are very distant, around $300\;km$ in all transformations, while the sources of Drava, Isel and Tauernbach are all much closer to the source of transformed Ister in all transformations, e.g. in Table \ref{tab:sources_JC} all distances are in the range from around 20 km to around 45 km, and they are slightly bigger in other Tables.

Now, let us look to Tables 3.x.2. Both the maximal and the mean Hausdorff distances (HD) are bigger when comparing Danube and Ister than when comparing Ister to the other river courses in their full length, see the first part of the Tables (first three rows). This difference is significantly emphasized in the second part of the Tables (fourth to the sixth row), where only the partial upper river courses are compared. The Hausdorff distances of Danube river and the reconstructed Ister on the upper part of their courses (HD of D1 and I1) are really high - the maximal HD is about 300 km and the mean HD is about 150 km. Opposite to that fact, the maximal and the mean HD of the transformed Ister and Drava/Tauernbach-Isel-Drava(TID) course are much lower. For example, in Table \ref{tab:HD_JC} the maximal Hausdorff distances are about 40 km and the mean Hausdorff distances are about 20 km only, which quantitatively express the visual similarity of the transformed Ister and Drava/TID courses in Figure \ref{fig:LAGL_JC}.

Tables 3.x.3 show the matching lengths in three different narrow bands 10, 50 and 100 km. It is another way to show how closely are the river streams on their length. The longer common length in a narrow band the better correspondence of the river courses is detected. As we see again in the second parts of these Tables (fifth to the sixth row), the matching length of upper courses of the transformed Ister and Drava/TID rivers is very high in narrow band 100 km (close or equal to 100\%) for all three LAGL transformations. The matching length is also very high in 50 km narrow band for the first two experiments and there is a similarity of river courses also in 10 km narrow band in the first experiment, which again show the perfect correspondence of the upper Ister and Drava/TID courses. On the other hand, there is almost no similarity of the transformed Ister and the Danube river in its upper course as seen in the fourth row of the Tables. 

The third part (the seventh row) of Tables 3.x.2 and 3.x.3 evaluate the quality of LAGL transformation of the Ister river. Since there is no doubt that the lower course of Danube and the transformed Ister should correspond to each other, the transformation which gives the lowest Hausdorff distances and the highest matching length on these partial river courses is the most reliable concerning the accuracy of the Strabo's Ister reconstruction. As one can see, from this point of view the most accurate is the LAGL transformation from Experiment 3.1, using just the Adriatic and the Black Sea coast regions. From the above discussion, we see that it also gives the best Strabo's Ister and Tauerbnbach-Isel-Drava-Danube correspondence. 
%Such result is interesting also from the mathematical modelling point of view, it uses the least possible information for LAGL, just two polygons on two respective sides of the Strabo's Ister, so the model is not over-fitted neither under-fitted, which may happen if we use more LAGL polygons or just one global affine transformation. 

\startsubsectionnumbering % ine cislovanie tabuliek a obrazkov
%\newpage
\invisiblesubsection{Transformation using the Adriatic and the Black Sea coast regions}

\clearpage

\begin{table}[h!] 
 \caption{The distance of river sources.\label{tab:sources_JC}} 
 \begin{center} \footnotesize 
 \begin{tabular}{|c|c|c|c|}\hline Danube - Ister & Drava - Ister & Isel - Ister & Tauernbach - Ister \\\hline \hline 
317.781 km & 29.849 km & 21.776 km &45.109 km \\ 
\hline \end{tabular} 
 \end{center} 
\end{table} 

\begin{table}[h!] 
 \caption{The table contains three parts. In the first part, the Hausdorff distances (HD) of the river courses in their full length are presented. Then the partial river courses are compared.\label{tab:HD_JC}} 
 \begin{center} \footnotesize 
 \begin{tabular}{|cc|d{3.3}c|d{3.3}c|} \hline \multicolumn{2}{|c|} {Curve} &\multicolumn{2}{c|} {Maximal HD[$km$]} &\multicolumn{2}{c|} {Mean HD[$km$]}\\	$A$& $B$& \multicolumn{1}{c}{${\delta}_{H}(A,B)$} & $ { \text{ d } }_{ H }(A, B)$& \multicolumn{1}{c}{$\overline{ \delta }_{ H }(A, B)$} & $\overline{ \text{d} }_{ H }(A, B)$\\\hline\hline
D & Ister & 317.781 & \multirow{2}{*}{317.781} & 112.802 & \multirow{2}{*}{100.898}\\
Ister & D & 206.829 & & 79.573& \\\hline
DD & Ister & 154.611 & \multirow{2}{*}{154.611} & 44.903 & \multirow{2}{*}{$\,\ $42.686}\\
Ister & DD & 149.693 & & 39.621& \\\hline
TIDD & Ister & 154.611 & \multirow{2}{*}{154.611} & 44.717 & \multirow{2}{*}{$\,\ $42.524}\\
Ister & TIDD & 149.693 & & 39.484& \\\hline
\hline D1 & I1 & 317.781 & \multirow{2}{*}{317.781} & 177.243 & \multirow{2}{*}{163.057}\\
I1 & D1 & 206.829 & & 130.218& \\\hline
Drava & I1 & 39.020 & \multirow{2}{*}{$\,\ $39.020} & 21.831 & \multirow{2}{*}{$\,\ $21.246}\\
I1 & Drava & 38.809 & & 20.547& \\\hline
TID & I1 & 39.020 & \multirow{2}{*}{$\,\ $40.628} & 21.383 & \multirow{2}{*}{$\,\ $20.833}\\
I1 & TID & 40.628 & & 20.173& \\\hline
\hline D2 & I2 & 154.611 & \multirow{2}{*}{154.611} & 65.344 & \multirow{2}{*}{$\,\ $62.953}\\
I2 & D2 & 149.693 & & 59.130& \\\hline
\end{tabular} 
 \end{center} 
\end{table} 

\begin{table}[h!] 
 \caption{The table contains three parts. In the first part, the matching lengths of the river courses in their full length are presented. Then the partial river courses are compared.\label{tab:ML_JC}} 
 \begin{center} \footnotesize 
 \scalemath{0.99}{
\begin{tabular}{|ccc|d{4.3}d{2.4}|d{4.3}d{2.4}|d{4.3}d{2.4}|} \hline \multicolumn{3}{|c|}{Curve} &\multicolumn{6}{c|}{Matching length $L_m(A, B, d_t) \; \left[km \right]$} \\\textit{A} & $L_A[km]$ & \textit{B} & \multicolumn{2}{c}{$d_t = 10$} & \multicolumn{2}{c}{$d_t = 50$} & \multicolumn{2}{c|}{$d_t = 100$} \\ \hline \hline D & 2714.887 & Ister & 228.127& (8.4$\%$) & 633.263& (23.3$\%$) & 1114.154& (41.0$\%$) \\ 
Ister & 1421.117 & D & 185.952 & (13.0$\%$) & 493.414 & (34.7$\%$) & 782.470 & (55.0$\%$) \\ 
\multicolumn{3}{| c |}{ A v e r a g e } & 207.040 & (10.0$\%$) & 563.339 & (27.2$\%$) & 948.312 & (45.8$\%$) \\ \hline 
DD & 2021.995 & Ister & 299.190& (14.7$\%$) & 1285.049& (63.5$\%$) & 1693.837& (83.7$\%$) \\ 
Ister & 1421.117 & DD & 267.783 & (18.8$\%$) & 1011.383 & (71.1$\%$) & 1247.252 & (87.7$\%$) \\ 
\multicolumn{3}{| c |}{ A v e r a g e } & 283.486 & (16.4$\%$) & 1148.216 & (66.6$\%$) & 1470.545 & (85.4$\%$) \\ \hline 
TIDD & 2026.547 & Ister & 320.930& (15.8$\%$) & 1289.600& (63.6$\%$) & 1698.389& (83.8$\%$) \\ 
Ister & 1421.117 & TIDD & 287.667 & (20.2$\%$) & 1011.383 & (71.1$\%$) & 1247.252 & (87.7$\%$) \\ 
\multicolumn{3}{| c |}{ A v e r a g e } & 304.299 & (17.6$\%$) & 1150.492 & (66.7$\%$) & 1472.821 & (85.4$\%$) \\ \hline 
\hline D1 & 1430.501 & I1 & 16.857& (1.1$\%$) & 85.824& (5.9$\%$) & 157.926& (11.0$\%$) \\ 
I1 & $\,\ $617.966 & D1 & 24.146 & (3.9$\%$) & 99.998 & (16.1$\%$) & 153.185 & (24.7$\%$) \\ 
\multicolumn{3}{| c |}{ A v e r a g e } & 20.501 & (2.0$\%$) & 92.911 & (9.0$\%$) & 155.555 & (15.1$\%$) \\ \hline 
Drava & $\,\ $737.609 & I1 & 87.920& (11.9$\%$) & 737.609& (100.0$\%$) & 737.609& (100.0$\%$) \\ 
I1 & $\,\ $617.966 & Drava & 105.976 & (17.1$\%$) & 617.966 & (100.0$\%$) & 617.966 & (100.0$\%$) \\ 
\multicolumn{3}{| c |}{ A v e r a g e } & 96.948 & (14.3$\%$) & 677.788 & (100.0$\%$) & 677.788 & (100.0$\%$) \\ \hline 
TID & $\,\ $742.161 & I1 & 109.660& (14.7$\%$) & 742.161& (100.0$\%$) & 742.161& (100.0$\%$) \\ 
I1 & $\,\ $617.966 & TID & 125.860 & (20.3$\%$) & 617.966 & (100.0$\%$) & 617.966 & (100.0$\%$) \\ 
\multicolumn{3}{| c |}{ A v e r a g e } & 117.760 & (17.3$\%$) & 680.063 & (100.0$\%$) & 680.063 & (100.0$\%$) \\ \hline 
\hline D2 & 1284.386 & I2 & 211.270& (16.4$\%$) & 547.439& (42.6$\%$) & 956.228& (74.4$\%$) \\ 
I2 & $\,\ $803.151 & D2 & 161.806 & (20.1$\%$) & 393.416 & (48.9$\%$) & 629.285 & (78.3$\%$) \\ 
\multicolumn{3}{| c |}{ A v e r a g e } & 186.538 & (17.8$\%$) & 470.428 & (45.0$\%$) & 792.757 & (75.9$\%$) \\ \hline 
\end{tabular} }
 \end{center} 
\end{table}

\clearpage
\invisiblesubsection{Transformation using the Adriatic coast, the Black Sea coast and the Greece and Albania regions}

\begin{table}[h!] 
 \caption{The distance of river sources.\label{tab:sources_JGC}} 
 \begin{center} \footnotesize 
 \begin{tabular}{|c|c|c|c|}\hline Danube - Ister & Drava - Ister & Isel - Ister & Tauernbach - Ister \\\hline \hline 
278.962 km & 67.198 km & 51.800 km &72.642 km \\ 
\hline \end{tabular} 
 \end{center} 
\end{table}

\begin{table}[h!] 
 \caption{The table contains three parts. In the first part, the Hausdorff distances (HD) of the river courses in their full length are presented. Then the partial river courses are compared.} 
 \begin{center} \footnotesize 
 \begin{tabular}{|cc|d{3.3}c|d{3.3}c|} \hline \multicolumn{2}{|c|} {Curve(river)} &\multicolumn{2}{c|} {Maximal HD[$km$]} &\multicolumn{2}{c|} {Mean HD[$km$]}\\	$A$& $B$& \multicolumn{1}{c}{${\delta}_{H}(A,B)$} & $ { \text{ d } }_{ H }(A, B)$& \multicolumn{1}{c}{$\overline{ \delta }_{ H }(A, B)$} & $\overline{ \text{d} }_{ H }(A, B)$\\\hline\hline
D & Ister & 278.962 & \multirow{2}{*}{278.962} & 116.695 & \multirow{2}{*}{108.088}\\
Ister & D & 194.446 & & 92.584& \\\hline
DD & Ister & 198.848 & \multirow{2}{*}{198.848} & 66.131 & \multirow{2}{*}{$\,\ $62.729}\\
Ister & DD & 189.915 & & 57.998& \\\hline
TIDD & Ister & 198.848 & \multirow{2}{*}{198.848} & 65.625 & \multirow{2}{*}{$\,\ $61.889}\\
Ister & TIDD & 189.915 & & 56.683& \\\hline
\hline D1 & I1 & 278.962 & \multirow{2}{*}{278.962} & 157.633 & \multirow{2}{*}{147.751}\\
I1 & D1 & 194.446 & & 125.898& \\\hline
Drava & I1 & 54.656 & \multirow{2}{*}{$\,\ $66.735} & 35.270 & \multirow{2}{*}{$\,\ $35.435}\\
I1 & Drava & 66.735 & & 35.623& \\\hline
TID & I1 & 54.656 & \multirow{2}{*}{$\,\ $70.863} & 33.857 & \multirow{2}{*}{$\,\ $33.083}\\
I1 & TID & 70.863 & & 32.194& \\\hline
\hline D2 & I2 & 198.848 & \multirow{2}{*}{198.848} & 96.661 & \multirow{2}{*}{$\,\ $93.564}\\
I2 & D2 & 189.915 & & 88.555& \\\hline
\end{tabular} 
 \end{center} 
\end{table} 

\begin{table}[h!] 
 \caption{The table contains three parts. In the first part, the matching lengths of the river courses in their full length are presented. Then the partial river courses are compared.} 
 \begin{center} \footnotesize 
 \scalemath{0.99}{
\begin{tabular}{|ccc|d{4.3}d{2.4}|d{4.3}d{2.4}|d{4.3}d{2.4}|} \hline \multicolumn{3}{|c|}{Curve(river)} &\multicolumn{6}{c|}{Matching length $L_m(A, B, d_t) \; \left[km \right]$} \\ \textit{A} & $L_A[km]$ & \textit{B} & \multicolumn{2}{c}{$d_t = 10$} & \multicolumn{2}{c}{$d_t = 50$} & \multicolumn{2}{c|}{$d_t = 100$} \\ \hline \hline D & 2714.887 & Ister & 61.260& (2.2$\%$) & 481.784& (17.7$\%$) & 948.962& (34.9$\%$) \\ 
Ister & 1441.085 & D & 55.447 & (3.8$\%$) & 323.687 & (22.4$\%$) & 641.770 & (44.5$\%$) \\ 
\multicolumn{3}{| c |}{ A v e r a g e } & 58.353 & (2.8$\%$) & 402.736 & (19.3$\%$) & 795.366 & (38.2$\%$) \\ \hline 
DD & 2021.995 & Ister & 26.097& (1.2$\%$) & 1082.523& (53.5$\%$) & 1508.358& (74.5$\%$) \\ 
Ister & 1441.085 & DD & 26.768 & (1.8$\%$) & 852.692 & (59.1$\%$) & 1144.834 & (79.4$\%$) \\ 
\multicolumn{3}{| c |}{ A v e r a g e } & 26.433 & (1.5$\%$) & 967.608 & (55.8$\%$) & 1326.596 & (76.6$\%$) \\ \hline 
TIDD & 2026.547 & Ister & 41.406& (2.0$\%$) & 1087.075& (53.6$\%$) & 1512.910& (74.6$\%$) \\ 
Ister & 1441.085 & TIDD & 50.989 & (3.5$\%$) & 855.911 & (59.3$\%$) & 1144.834 & (79.4$\%$) \\ 
\multicolumn{3}{| c |}{ A v e r a g e } & 46.197 & (2.6$\%$) & 971.493 & (56.0$\%$) & 1328.872 & (76.6$\%$) \\ \hline 
\hline D1 & 1430.501 & I1 & 34.437& (2.4$\%$) & 112.681& (7.8$\%$) & 178.212& (12.4$\%$) \\ 
I1 & $\,\ $646.869 & D1 & 25.465 & (3.9$\%$) & 90.754 & (14.0$\%$) & 143.805 & (22.2$\%$) \\ 
\multicolumn{3}{| c |}{ A v e r a g e } & 29.951 & (2.8$\%$) & 101.717 & (9.7$\%$) & 161.008 & (15.5$\%$) \\ \hline 
Drava & $\,\ $737.609 & I1 & 0.000& (0.0$\%$) & 713.419& (96.7$\%$) & 737.609& (100.0$\%$) \\ 
I1 & $\,\ $646.869 & Drava & 0.000 & (0.0$\%$) & 619.759 & (95.8$\%$) & 646.869 & (100.0$\%$) \\ 
\multicolumn{3}{| c |}{ A v e r a g e } & 0.000 & (0.0$\%$) & 666.589 & (96.2$\%$) & 692.239 & (100.0$\%$) \\ \hline 
TID & $\,\ $742.161 & I1 & 15.308& (2.0$\%$) & 717.971& (96.7$\%$) & 742.161& (100.0$\%$) \\ 
I1 & $\,\ $646.869 & TID & 24.221 & (3.7$\%$) & 622.978 & (96.3$\%$) & 646.869 & (100.0$\%$) \\ 
\multicolumn{3}{| c |}{ A v e r a g e } & 19.764 & (2.8$\%$) & 670.474 & (96.5$\%$) & 694.515 & (100.0$\%$) \\ \hline 
\hline D2 & 1284.386 & I2 & 26.097& (2.0$\%$) & 369.103& (28.7$\%$) & 770.749& (60.0$\%$) \\ 
I2 & $\,\ $794.215 & D2 & 26.768 & (3.3$\%$) & 232.933 & (29.3$\%$) & 497.964 & (62.6$\%$) \\ 
\multicolumn{3}{| c |}{ A v e r a g e } & 26.433 & (2.5$\%$) & 301.018 & (28.9$\%$) & 634.357 & (61.0$\%$) \\ \hline 
\end{tabular} }
 \end{center} 
\end{table}

\clearpage
\invisiblesubsection{Transformation using the Adriatic coast, the Black sea coast and the Alps and Italy regions}
 
\begin{table}[h!] 
 \caption{The distance of river sources.\label{tab:sources_AJC}} 
 \begin{center} \footnotesize 
 \begin{tabular}{|c|c|c|c|}\hline Danube - Ister & Drava - Ister & Isel - Ister & Tauernbach - Ister \\\hline \hline 
296.586 km & 50.094 km & 35.885 km &58.252 km \\ 
\hline \end{tabular} 
 \end{center} 
\end{table}

\begin{table}[h!] 
 \caption{The table contains three parts. In the first part, the Hausdorff distances (HD) of the river courses in their full length are presented. Then the partial river courses are compared.} 
 \begin{center} \footnotesize 
 \begin{tabular}{|cc|d{3.3}c|d{3.3}c|} \hline \multicolumn{2}{|c|} {Curve(river)} &\multicolumn{2}{c|} {Maximal HD[$km$]} &\multicolumn{2}{c|} {Mean HD[$km$]}\\	$A$& $B$& \multicolumn{1}{c}{${\delta}_{H}(A,B)$} & $ { \text{ d } }_{ H }(A, B)$& \multicolumn{1}{c}{$\overline{ \delta }_{ H }(A, B)$} & $\overline{ \text{d} }_{ H }(A, B)$\\\hline\hline
D & Ister & 296.586 & \multirow{2}{*}{296.586} & 117.455 & \multirow{2}{*}{108.467}\\
Ister & D & 205.818 & & 92.198& \\\hline
DD & Ister & 219.255 & \multirow{2}{*}{219.255} & 84.164 & \multirow{2}{*}{$\,\ $80.301}\\
Ister & DD & 205.818 & & 74.903& \\\hline
TIDD & Ister & 219.255 & \multirow{2}{*}{219.255} & 83.610 & \multirow{2}{*}{$\,\ $79.679}\\
Ister & TIDD & 205.818 & & 74.175& \\\hline
\hline D1 & I1 & 296.586 & \multirow{2}{*}{296.586} & 145.940 & \multirow{2}{*}{137.085}\\
I1 & D1 & 194.015 & & 115.984& \\\hline
Drava & I1 & 101.197 & \multirow{2}{*}{101.197} & 65.823 & \multirow{2}{*}{$\,\ $63.574}\\
I1 & Drava & 94.308 & & 60.811& \\\hline
TID & I1 & 101.197 & \multirow{2}{*}{101.197} & 64.181 & \multirow{2}{*}{$\,\ $61.763}\\
I1 & TID & 94.308 & & 58.773& \\\hline
\hline D2 & I2 & 219.255 & \multirow{2}{*}{219.255} & 112.879 & \multirow{2}{*}{107.675}\\
I2 & D2 & 205.818 & & 99.975& \\\hline
\end{tabular} 
 \end{center} 
\end{table} 

\begin{table}[h!] 
 \caption{The table contains three parts. In the first part, the matching lengths of the river courses in their full length are presented. Then the partial river courses are compared.} 
 \begin{center} \footnotesize 
 \scalemath{0.99}{
\begin{tabular}{|ccc|d{4.3}d{2.4}|d{4.3}d{2.4}|d{4.3}d{2.4}|} \hline \multicolumn{3}{|c|}{Curve(river)} &\multicolumn{6}{c|}{Matching length $L_m(A, B, d_t) \; \left[km \right]$} \\\textit{A} & $L_A[km]$ & \textit{B} & \multicolumn{2}{c}{$d_t = 10$} & \multicolumn{2}{c}{$d_t = 50$} & \multicolumn{2}{c|}{$d_t = 100$} \\ \hline \hline D & 2714.887 & Ister & 46.134& (1.6$\%$) & 327.386& (12.0$\%$) & 911.606& (33.5$\%$) \\ 
Ister & 1468.316 & D & 40.573 & (2.7$\%$) & 248.261 & (16.9$\%$) & 674.137 & (45.9$\%$) \\ 
\multicolumn{3}{| c |}{ A v e r a g e } & 43.354 & (2.0$\%$) & 287.824 & (13.7$\%$) & 792.871 & (37.9$\%$) \\ \hline 
DD & 2021.995 & Ister & 22.798& (1.1$\%$) & 374.882& (18.5$\%$) & 1390.515& (68.7$\%$) \\ 
Ister & 1468.316 & DD & 20.561 & (1.4$\%$) & 301.285 & (20.5$\%$) & 1130.088 & (76.9$\%$) \\ 
\multicolumn{3}{| c |}{ A v e r a g e } & 21.679 & (1.2$\%$) & 338.083 & (19.3$\%$) & 1260.301 & (72.2$\%$) \\ \hline 
TIDD & 2026.547 & Ister & 48.427& (2.3$\%$) & 379.434& (18.7$\%$) & 1395.067& (68.8$\%$) \\ 
Ister & 1468.316 & TIDD & 48.975 & (3.3$\%$) & 296.829 & (20.2$\%$) & 1130.088 & (76.9$\%$) \\ 
\multicolumn{3}{| c |}{ A v e r a g e } & 48.701 & (2.7$\%$) & 338.131 & (19.3$\%$) & 1262.577 & (72.2$\%$) \\ \hline 
\hline D1 & 1430.501 & I1 & 19.581& (1.3$\%$) & 122.200& (8.5$\%$) & 257.106& (17.9$\%$) \\ 
I1 & $\,\ $600.305 & D1 & 10.404 & (1.7$\%$) & 57.143 & (9.5$\%$) & 144.355 & (24.0$\%$) \\ 
\multicolumn{3}{| c |}{ A v e r a g e } & 14.992 & (1.4$\%$) & 89.671 & (8.8$\%$) & 200.730 & (19.7$\%$) \\ \hline 
Drava & $\,\ $737.609 & I1 & 0.000& (0.0$\%$) & 191.060& (25.9$\%$) & 727.479& (98.6$\%$) \\ 
I1 & $\,\ $600.305 & Drava & 0.000 & (0.0$\%$) & 178.483 & (29.7$\%$) & 600.305 & (100.0$\%$) \\ 
\multicolumn{3}{| c |}{ A v e r a g e } & 0.000 & (0.0$\%$) & 184.771 & (27.6$\%$) & 663.892 & (99.2$\%$) \\ \hline 
TID & $\,\ $742.161 & I1 & 25.628& (3.4$\%$) & 195.612& (26.3$\%$) & 732.031& (98.6$\%$) \\ 
I1 & $\,\ $600.305 & TID & 28.413 & (4.7$\%$) & 174.027 & (28.9$\%$) & 600.305 & (100.0$\%$) \\ 
\multicolumn{3}{| c |}{ A v e r a g e } & 27.021 & (4.0$\%$) & 184.819 & (27.5$\%$) & 666.168 & (99.2$\%$) \\ \hline 
\hline D2 & 1284.386 & I2 & 22.798& (1.7$\%$) & 183.821& (14.3$\%$) & 652.906& (50.8$\%$) \\ 
I2 & $\,\ $868.010 & D2 & 20.561 & (2.3$\%$) & 122.801 & (14.1$\%$) & 529.782 & (61.0$\%$) \\ 
\multicolumn{3}{| c |}{ A v e r a g e } & 21.679 & (2.0$\%$) & 153.311 & (14.2$\%$) & 591.344 & (54.9$\%$) \\ \hline 
\end{tabular} }
 \end{center} 
\end{table} 

%\clearpage

\startsectionnumbering
\section{Discussion on some historical issues} We could finish our work here, just by developing the mathematical model and the numerical algorithm and by showing the correspondence of the Strabo's Ister with the nowadays Tauernbach-Isel-Drava-Danube course. But, since any new result in the mathematical and computational modelling used to bring new insight into the related pure or applied science problem, we are going to do the same for the history which forms our application background.

As we have already stated in the Introduction, many historical claims assuming that in the time of Strabo the upper course of the river Ister corresponds to the nowadays upper course of the Danube river must be revisited. Regarding this fact, there are many interesting open questions arising but two of them are the most interesting for us, {\it (i)}
where was then located the so-called Hercynian Forest/Hercynia silva/Herkynian Forest (Ἑρκυνίου δρυμος)
\cite{Lacus, MD, Roller}, the seat of Suevi/Suevi/Soebians (Σοήβων) \cite{Lacus, MD, Roller} ? And {\it (ii)} who were Strabo's {\bf Suevi}? 

In the sequel, we will use the terms Suevi and Hercynian Forest because they seem to be the most spread in the English and Latin literature. We also note that in the Slovak literature the terms like Sv\'ebi/Suavi are used as well \cite{Steinhubel}. First of all, at the end of section (4-6-9) Strabo writes that the Ister source is near the seats of Suevi and the Hercynian Forest: 

"ὅπου αἱ τοῦ Ἴστρου πηγαὶ πλησίον Σοήβων καὶ τοῦ Ἑρκυνίου δρυμοῦ". 

Since we have shown above where the Ister source is located by means of Strabo's {\it Geographica}, we can clearly state that Strabo's Suevi near the Ister source have no relation with the Swabia (and the Swabians) in Bavaria, Germany, but we can claim that Strabo is speaking here about a settlement in the south-east Alpine region, around the boundaries of nowadays Carinthia, Tyrol, north-east Italy and Slovenia. With a high probability, this settlement was Slavic in that time and before which can be confirmed by many geographic names around the Strabo's Ister source which are of Slavic origin. 
The origin of local geographic names in the neighbourhood of Val di Pusteria, in the valleys of upper Drava, Villgraten, Gail and Isel rivers, was studied in \cite{Alpy} where almost 200 names from this local area, including settlements, rivers and creeks, hills, forests or meadows, of the Slavic origin were presented. This study is based on works \cite{Miklosic1, Miklosic2} by Franc Miklo\v si\v c (Franz Miklosich), one of the most respected philologists of the Habsburg empire in the second half of the 19th century. In two volumes Miklo\v si\v c presented all the important rules for creating the Slavic geographic names (Vol. I) and he collected a comprehensive set of 789 bases of Slavic geographic names (Vol. II) from the whole Habsburg empire. He also gave the most common rules for changing Slavic names to German (and Hungarian) such as change of the Slavic "B" to German "F", etc. Now, we present just a few examples of geographic names of the Slavic origin around the Ister source. Here, and also in further paragraphs of this section, we give the meaning of these geographic names together with English also in the Slovak language, because it is the most familiar to authors and there does not exist any common Slavic language, and where it has a sense we give it also in the Slovenian due to \cite{Alpy}.

First interesting geographic name is Val di Pusteria (Pustertal in German, Puster Valley in English), which would be in Slovak "Pust\'e \'udolie" or "Pust\'a dolina" and in Slovenian "Pusta dolina" or "Pustodol", see \cite{Alpy} and \cite{Miklosic2} - point 512, meaning "Deserted Valley" in English. At the eastern end of the western (lower) part of the Puster Valley, there is the castle hill, nowadays called Heinfels, below which the river Villgraten(bach) empties to the Drava. By \cite{Alpy} the river name has the Slavic origin "Velegrad" and there are two other creeks around, Gradenbach and Gratzbach (mentioned as Gradiz in \cite{Alpy}) with the same Slavic base "grad", see \cite{Miklosic2}-122. A bit to the west, there is the village Versciaco (Vierschach in German) with the meaning "V\'r\v sok" in Slovak and "hillock" in English, see also \cite{Alpy}.
Going further north-east, in the Defereggen Valley there is the settlement Feistritz, nowadays part of St. Jacob in Defereggen, and also the creek Feistritz(bach). Feistritz is a German writing of the Slavic name "Bystrica" in Slovak and "Bistrica" in Slovenian, see \cite{Alpy}, \cite{Miklosic2}-45. As a further nice example we mention Pro{\ss}egg(-klamm) \cite{Alpy}, the village and the gorge on the Tauernbach creek north of the town Matrei in Osttirol. The name Pro{\ss}egg has exact analogy in Prosiek village and gorge in Cho\v csk\'e vrchy, Slovakia, since "priesek" in Slovak just means "gorge" in English and "klamm" in German. Another interesting fact is given by the historical names of Matrei in Osttirol town and Gro{\ss}venediger mountain. They were called Windisch Matray als Mauter and Windisch Taurn, respectively, on the map in Figure \ref{fig:isola}, and it is generally accepted that "Windisch" meant Slavic in German dialects. Further north, there is the valley and the long creek Frosnitz(bach) \cite{Alpy}, with the source just below the Gro{\ss}venediger glaciers, emptying to Tauernbach near Gruben village. The name Frosnitz is a German writing of the Slavic name "Brusnica" meaning "cranberry" in English, see \cite{Alpy} and \cite{Miklosic2}-33. Interestingly, another creek with the same name meaning, the Fruschnitz(bach), is stemming from the western side of Gro{\ss}glockner \cite{Alpy}.
%On the south, above the Val di Pusteria, there is the Rocca dei Baranci, and "baranci" means group of Aries in the Slovak language, see also \cite{Miklosic2}-11, thus the name in Slovak would be "Skala barancov" or "Baran\v cia skala" or just "Baranci", meaning in English the Aries' rock - expressing the shape of the rocks in the form of Aries horns. Horn means "roh" in the Slovak language and there is an analogy in the name of Roh\'a\v ce mountains containing Baranec hill in the West Tatra Mountains and there is also the mountain Baranie Rohy in the High Tatra mountains. We are aware of the fact that "barancio" also means {\it Pinus mugo} in Italian, but the Slavic meaning of the name is convincing, too, thus we present it. 
%A further example is the name of Drava river itself, with a clear meaning and almost exactly the same writing in Slovak, as "Drav\'a rieka" meaning the "Ravenous river" (in the sense of flow). 
%It could be the name of a local Ister (Tauernbach-Isel-Drava-Danube) tributary represented by the Drava course until the confluence with Isel because Strabo mentions in section (4-6-9) the tributary with a similar name Δούρας. Another interesting option (daring hypothesis) comes from the Strabo's Greek form of the Ister name - Istros (Ἴστρος) which could be also derived from "Bystr\'y - Bystros\v t" (quickly flowing) property with a very similar meaning as "Drav\'y - Dravos\v t" (ravenous) property, see also the next paragraphs of this section. 
But the most astonishing is the name of the longest glacier in the Eastern Alpine region - Pasterze - in Slovak pronunciation "Pastierce" which means the place of shepherds or the pastureland. Such form of the place name has the classical Slavic suffix -ce (-ze in German writing), see \cite{Miklosic1} - Chapter 2, Section V, points 18 and 17. 
%and many analogies in Slavic countries, e.g. in Slovakia, just to mention few - V\v celince, Kovarce, Lovce, Zlievce, Plachtince, \v Ce\v ladince, Hos\v tovce, Pl\'a\v s\v tovce, Hudcovce, Sedlice, Hrn\v ciarovce, Mlyn\'arovce, Medovarce - always expressing some specific "homogeneous" activity performed by "craftsmen" or another group of inhabitants on that place or in a village. 
And why is it so astonishing? Recent results studying peat samples from the area of the retreating Pasterze glacier indicates grass-like pastureland vegetation (Cyperaceae-{\it Carex}, {\it Bidens alba} - shepherd's needles) and human impact on the vegetation during the Subboreal  Chronozone (3780-800 BC) a warmer period  of the Holocene \cite{Pasterze}. In that period only, the name Pasterze with the Slavic pastureland meaning could be given to that place by the inhabitants of this Alpine region - with a very high probability by the Suevi people. They lived there in the Strabo's times and there is no reason to assume that they did not settle there before. We adopt here the hypothesis of the continuity of settlement if no change is recorded in any historical source which is in agreement, e.g., with the Palaeolithic Continuity Theory of Mario Alinei declaring the stability and continuity of Romance-Celtic/Germanic/Slavic ethnic and language geographic distribution in Europe from the Upper Palaeolithic period \cite{Alinei, Alinei2}. It is also worth to note that such assumption, assuming stable Slavic sedentary population in this part of Europe already in the Strabo's times, does not exclude any immigration of the same or different ethnic origin and/or local acculturation in the later periods.

\begin{figure}[ht]
 \begin{center}
 \includegraphics[width = 120mm]{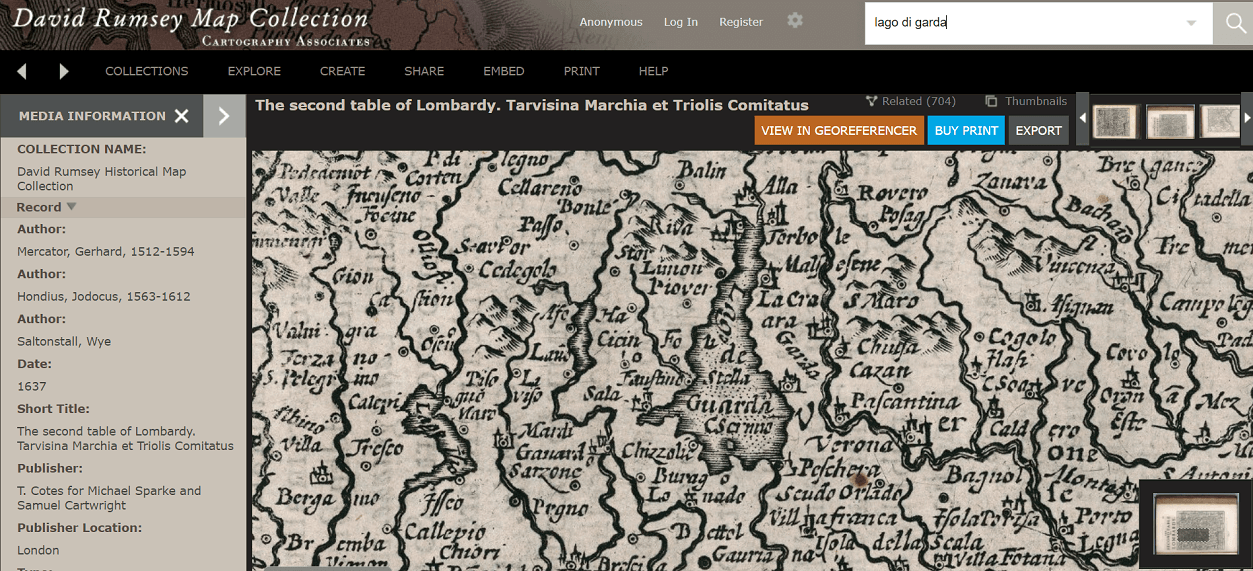}
 \caption{Detail of the map by Gerhard Mercator from 1639 \cite{Mercator2} where in the middle of the map one can see the Sirmione island in Lago di Garda, nowadays connected to the land, where Tiberius could build his military camp in 15 BC during the campaign against Rhaeti. 
}
 \label{fig:garda}
 \end{center}
\end{figure}

Next important mentions of the Suevi and Hercynian Forest is in section (7-1-5) of {\it Geographica} where the fights of Tiberius against the Rhaetians around 15 BC are mentioned, see also Cassius Dio \cite{Cassius} sections (54-22-1)-(54-22-5) and \cite{Legions}. For that military campaigns, Tiberius even built a military camp on the island - Sirmione - of the lake - Lago di Garda, see the map on Figure \ref{fig:garda} and \cite{Legions}. By Strabo, the lake is located south of the Ister source and close to the Rhaetian territories (which are around the Tridentine Alps by Cassius Dio \cite{Cassius}), fitting correctly the Lago di Garda and our location of the source of Ister. Then the journey to the Hercynian forest is described by: 

"ὥστ᾽ ἀνάγκη τῷ ἐκ τῆς Κελτικῆς ἐπὶ τὸν Ἑρκύνιον δρυμὸν ἰόντι πρῶτον μὲν διαπερᾶσαι τὴν λίμνην, ἔπειτα τὸν Ἴστρον, εἶτ᾽ ἤδη δι᾽ εὐπετεστέρων χωρίων ἐπὶ τὸν δρυμὸν τὰς προβάσεις ποιεῖσθαι δι᾽ ὀροπεδίων." 

So, from the (Cisalpine) Keltike, which is for Strabo the north part of nowadays Italy up to the base of the Alps (see section (5-1-3) of {\it Geographica}), one has to go along that lake (Lago di Garda) then continue up to and along the course of Ister (taking simply the route along the Adige-Isarco-Rienza-Drava courses) and then continue straightforwardly through more favourable upland planes to end up in the Hercynian Forest. Just looking at any map, e.g. Figure \ref{fig:LAGL_JC}, one clearly see that this journey must finish in the Carpathian-Alpine basin or better say in the large region between the Eastern Alps on the west and the Carpathian mountains ridge on the east and north. The last sentence of the section (7-1-5) is also very interesting: "ἔστι δὲ καὶ ἄλλη ὕλη μεγάλη Γαβρῆτα ἐπὶ τάδε τῶν Σοήβων, ἐπέκεινα δ᾽ ὁ Ἑρκύνιος δρυμός: ἔχεται δὲ κἀκεῖνος ὑπ᾽ αὐτῶν." It is not easy to translate, but in any case, it says that the whole Hercynian Forest, together with another large forest Gabreta (Γαβρῆτα), belongs to and is the seat of the Suevi. In the case that the Gabreta is on the west side of the Hercynian Forest it should correspond to the mountainous forested regions below the main Alpine ridges in the east and north-east directions such as the lower parts of Carinthia, Styria and nowadays Vienna Forest. In the case that the Gabreta is on the east and north sides
%, in other words beyond/on yonder side (ἐπέκεινα), 
of the Hercynian Forest, then it should correspond to the Carpathians mountains - Karpaty in the Slovak language. 
%Here we just notice a similarity of pronunciation of Gabreta and Karpaty because of the assimilation of K-G and b-p but do not favorize any hypothesis. 
A further indication of the correctness of the location of the Hercynian Forest in between the Alps and the Carpathians is given in section (7-3-1) of {\it Geographica} where Strabo mentions the land of Getians which were assumed to be a Thracians by Hellenes, see (7-3-2): "οἱ τοίνυν Ἕλληνες τοὺς Γέτας Θρᾷκας ὑπελάμβανον" and Thracia was generally considered as the land north of Greece up to the Ister (lower Danube). And in (7-3-1) Strabo explicitly says that the land of Getians extends along the southern side of the Istros, nowadays northern Bulgaria, and also on the opposite side, on the mountain slopes of the Hercynian Forest: 

"εἶτ᾽ εὐθὺς ἡ τῶν Γετῶν συνάπτει γῆ, κατ᾽ ἀρχὰς μὲν στενή, παρατεταμένη τῷ Ἴστρῳ κατὰ τὸ νότιον μέρος, κατὰ δὲ τοὐναντίον τῇ παρωρείᾳ τοῦ Ἑρκυνίου δρυμοῦ". 

From there, it is obvious that the Getians territory north of the Ister corresponds to the nowadays south part of Romania behind the Carpathian ridges and it is adjoining the Hercynian Forest. These facts indicate that the Hercynian Forest, the seat of Suevi, corresponds (at least) to the Carpathian-Alpine basin including the south-eastern Alps and Carpathian mountains as well. 

From the above facts, we can conclude that the Hercynian Forest - the Carpathian-Alpine basin in a broad sense - was in the times of Strabo settled by the Suevi people, the large ethnic group living in this compact area encompassed by the mountain ranges, as he says in sections (7-1-3) and (7-1-5) of {\it Geographica}. There are many remains of such compact settlement mainly in geographic names of rivers, mountains, towns and villages in the whole region, see e.g. \cite{Miklosic1,Miklosic2,Stanislav1, Stanislav2}. 
%First of all, just as a curiosity, we mention pronunciation similarity of "Hercynian/Herkynian" and "Uhorsko" (or "Uhersko" in the Moravian Slovak dialect and also in the nowadays Czech language), with a possible meaning in the Slovak language the land below (or close to) the mountains. The name "Uhorsko" has been used as a historical name for the Carpathian basin among Slovaks, see e.g. \cite{Sasinek}, and in a few other Slavic countries in some adjusted forms as well. Then, there are so many remaining common, exactly same or very similar in pronunciation and writing, geographic names around the whole Strabo's Hercynian Forest, see e.g. \cite{Miklosic1,Miklosic2,Stanislav1, Stanislav2}. 
For interested reader, we are going to touch some of them, concentrating our description mainly on the west (Slovenia, Carinthia, Styria), north (Slovakia) and east (Romania) nodes of an imaginary triangle in the Carpathian-Alpine basin.

Probably the most spread geographic name of the rivers and settlements are the names Bystrica/Bystr\'a (in Slovakia), Bistrica/Bistra (in Slovenia), Bistri{\textcommabelow{t}}a/ Bistra (in Romania) and Feistritz (in Austria), meaning "quickly flowing" in Slavic languages, see also \cite{Miklosic2}-45. In Slovakia, there are five Bystrica and one Bystr\'a settlements, e.g. towns Bansk\'a Bystrica, Pova\v zsk\'a Bystrica, etc., two river streams with the name Bystrica, one mountain peak Bystr\'a in the West Tatra mountains and there is also saddle Bystr\'e sedlo in the High Tatra mountains. In Slovenia, there are at least ten towns and villages with the name Bistrica, e.g. Ilirska Bistrica, Slovenska Bistrica, etc., the river Kamni\v ska Bistrica and stream Bistra. In Austria, there are at least eight settlements with the name Feistritz, mainly in Carinthia and Styria, there is a saddle Feistritz-Sattel, on the border of Styria and Lower Austria, below which is the spring of a long Styrian river Feistritz. We note that by \cite{Miklosic2,Alpy}, around 1870-1880 there was reported Feistritz 15 times in Carinthia and 40 times in Styria.  Surprisingly also in Romania, there are at least seven Bistri{\textcommabelow{t}}a towns and villages and five rivers with such name and even more, nine Bistra rivers and three such settlements, distributed all around the country, in counties Alba, Bacău, Bihor, Bistrița-Năsăud, Caraș-Severin, Gorj, Maramureș, Mehedinţi, Mureș, Neamț, Olt, Sibiu, Suceava, Vâlcea, and there is also a mountain range called Bistri{\textcommabelow{t}}a in northern central Romania. It is worth to note that there are many further such examples of common geographic names, e.g. "Slatina", with the meaning a "mineralized water" (\cite{Miklosic2}-585) and with the exact same writing at all the places, "Trnava" and its analogy "Târnava" in Romania, with the meaning of adjective related to "thorn" (\cite{Miklosic2}-696). Interestingly, in Romania, there are whole regions with very dense names of villages almost exclusively using the basis of the Slavic words or even more, copying almost exactly the village names used e.g. in Slovakia. For example, in Caraș-Severin county, near the Valea Cernei (\'Udolie \v Ciernej (rieky) in Slovak), are the villages Camena (Kamenn\'a in Slovak), Cozia (Kozia), Dobraia (Dobr\'a), Dolina (Dolina), Gruni (Gr\'u\v n), Hora Mare (Ve\v lk\'a Hora), Hora Mică (Mal\'a Hora), Iablanița (Jablonica), Ilova (Ilava), Obița (Obyce), Rusca (Rusk\'a), Ruștin (Hru\v st\'in), Sadova Nouă (Nov\'a Sadov\'a), Sadova Veche (Star\'a Sadov\'a),  Slatina-Timiș (Timi\v ssk\'a Slatina), Studena (Studen\'a), Topla (Tepl\'a), Zbegu (Zbehy), Plugova (Pluhov\'a), Zoina (Zoln\'a), etc., and there exist such examples around all Romania. 
Another interesting example is the usage of the geographic name "Studena" with the meaning "cold" (\cite{Miklosic2}-636) in practically the same form in Slovakia, Romania, Serbia, Croatia, Slovenia and even in Italy (Studena Alta and Studena Bassa villages in the province of Udine close to Carinthian border) although nowadays in the majority of these languages the word "studena" is not used to express the coldness. After this general overview, we touch further examples of Slavic geographic names which interest us because they appear in Roman and Greek writings from the beginning of the first millennium.

\begin{figure}[ht]
 \begin{center}
 \includegraphics[width = 1.\textwidth]{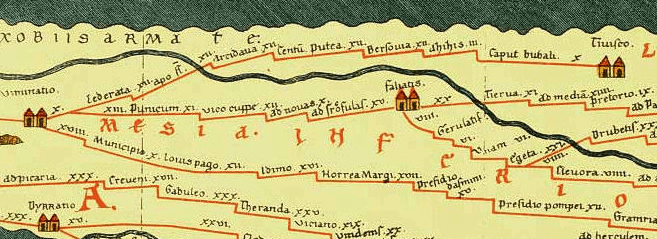} 
 \caption{A detail of Tabula Peutingeriana with the Trajan's roads to the province of Dacia and two stations "Bersouia" (on the top road) and "Tierua" (on the second from the top road). Source: Wikipedia.
}
 \label{fig:Peutingeriana1}
 \end{center}
\end{figure}

\begin{figure}[ht]
 \begin{center}
 \includegraphics[height = 60mm]{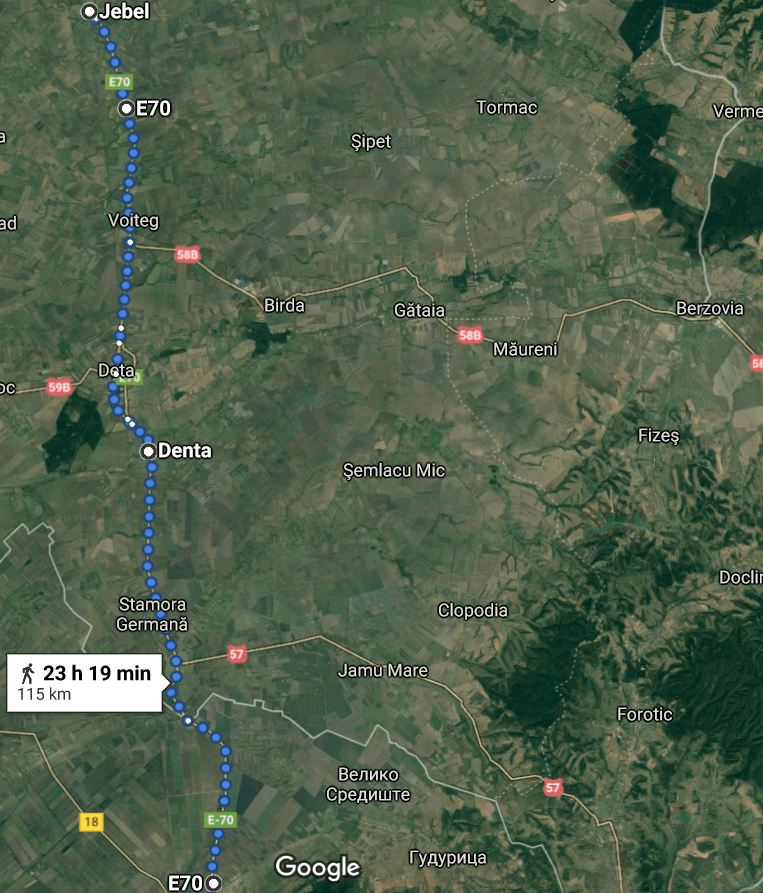}
 \hskip 3mm
 \includegraphics [height = 60mm]{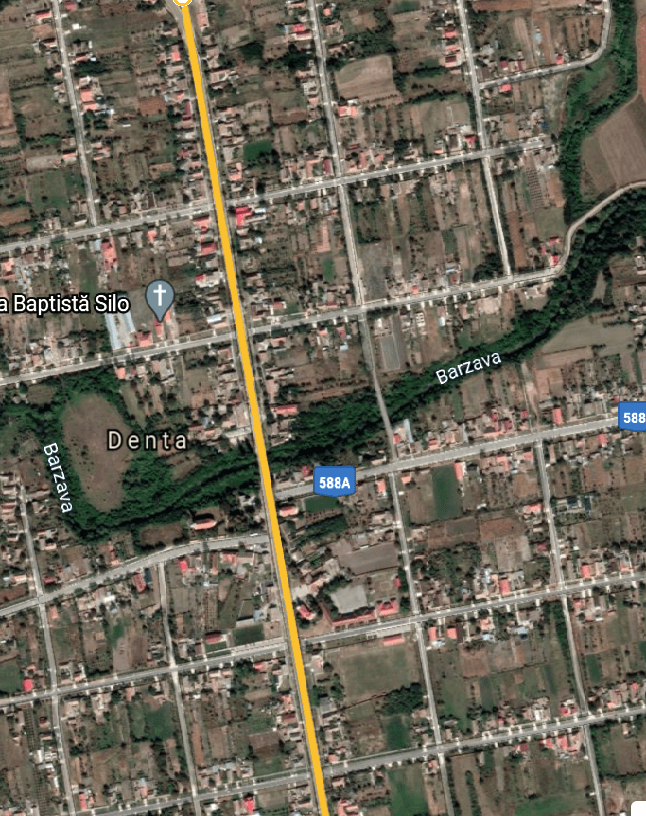}
 \caption{A detail of the estimated Trajan's road to "Tiuisco" (Timișoara) on Tabula Peutingeriana (left image) with four stops indicated, "Centu Putea" (E70 mark on the bottom), "Bersouia" (Denta), "Azizis" (E70 mark on the top) and "Caput bubali" (Jebel). Station "Bersovia" is on the road crossing with the Bârzava river in the village Denta (right image). The village Berzovia which is also located on the river Bârzava is more on the east and the road going through it would not fulfil the distances indicated on the Tabula Peutingeriana. That's why we think that the station called Berzovia was at the place where the road crossed the river with the same name. 
}
 \label{fig:Bersovia}
 \end{center}
\end{figure}

First, let us take the geographic names derived from the basis of Slavic word "breza" which means "birch" tree (\cite{Miklosic2}-29). In Slovakia, we have towns Brezov\'a pod Bradlom and Brezno and villages Brezovica (twice), Brezovi\v cka, Brezov, Rimavsk\'e Brezovo and \v Cesk\'e Brezovo and one river stream Brezovsk\'y potok (Brezovka). In Slovenia, there are villages Brezova, Brezova Reber pri Dvoru, Brezno (twice), Breza and Brezovo. In Romania there are several forms of this name, there are seven towns or villages and three rivers with the name Breaza, two villages and two river streams with the name Breazova, one of which is a tributary of the river Bârzava, which is another form of the same name. The river Bârzava, flowing through historical regions of Banat (Romania) and Vojvodina (Serbia), has further tributaries  Bârzăvița and Berzovița, there is the village Berzovia on the river Bârzava and very near is another village Brezon. It used to be claimed that the village Berzovia on the river Bârzava is noted on the Tabula Peutingeriana \cite{Tabula} as the station "Bersouia", see e.g. Pavol Jozef \v Saf\'arik seminal work \cite{Safarik}. Tabula Peutingeriana shows the Roman road system in the first centuries AD and it should be last revised in the 4th or early 5th century. The "Bersouia" is one of the stations on the most north Trajan's road to the province of Dacia which crosses the Ister (lower Danube) near the Banatska Palanka and continues in the direction to "Tiuisco" (Timișoara), see Fig. \ref{fig:Peutingeriana1}. In fact, the "Bersouia" is mentioned also in the Trajan's work {\it Dacica}, from around 100 AD, as "inde Berzobim, deinde Aizi processimus", meaning going from "Bersouia" to "Azizis", see Fig. \ref{fig:Peutingeriana1}. Taking into account the Tabula Peutingeriana distances given in Roman miles we claim that "Bersouia" or "Bersobis" is not directly the Berzovia village but it is (with high probability) the crossing point of that Trajan's road with the Bârzava river at a place near (or in) the Denta village, see Fig. \ref{fig:Bersovia}. Such claim, however, does not change anything on the interesting fact that the widely spread in Carpathian basin geographic name of the Slavic origin is mentioned on the map from the first centuries AD and in the work of emperor Trajan from around 100 AD.

%\begin{figure}[ht]
% \begin{center}
% \includegraphics[height = 40mm]{img/Bela_rieka.png} 
% \hskip 3mm
% \includegraphics [height = 40mm]{img/Fella_river.png}
% \caption{River Bel\'a in the High Tatras (left) and river Fella (Bela) in the Carnic Alps.
%}
% \label{fig:bela_fella}
% \end{center}
%\end{figure}

%\begin{figure}[ht]
% \begin{center}
% \includegraphics[width = 80mm]{img/Baile_Herculane_road_1824.png}
% \caption{Baile Herculane (Hercules baths) - near this famous spa town in the Caraș-Severin county, Romania, the river Belareka empties into the river Cerna. 
%}
% \label{fig:Belareka}
% \end{center}
%\end{figure}

\begin{figure}[ht]
 \begin{center}
 \includegraphics[width = 1.\textwidth]{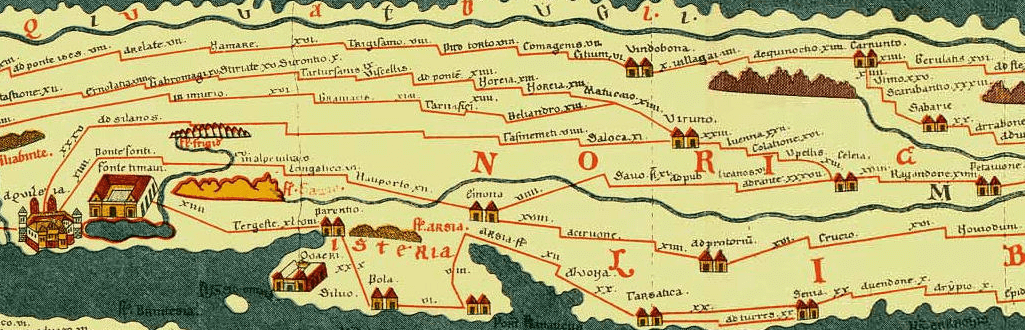}
 \caption{A detail of Tabula Peutingeriana with the road from "Petauione" through "Celeia" to "Beliandro" and farther away. Also stations "Tergeste" and "fonte timaui" can be seen on this map detail in the left bottom part. Source: Wikipedia.
}
 \label{fig:Peutingeriana2}
 \end{center}
\end{figure}

\begin{figure}[ht]
 \begin{center}
 \includegraphics[width = 1.\textwidth]{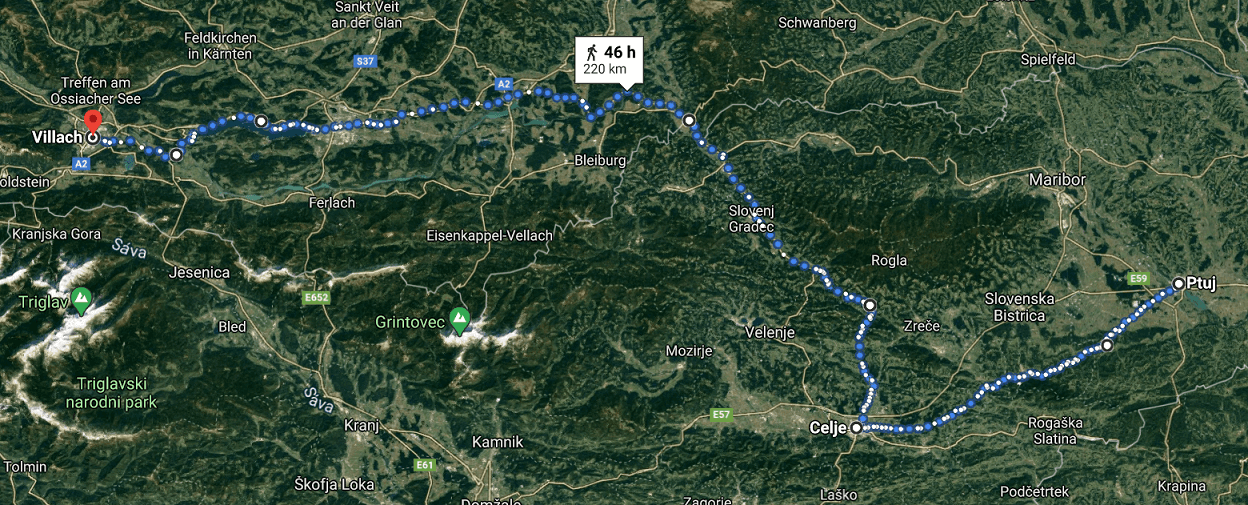}
 \caption{Distance from Ptuj to Villach through Celje is 220 km by the Google maps while by the Tabula Peutingeriana it should be 137 Roman miles which is approximately 222 km. Such accurate correspondence of distances shows that nowadays Villach (Beljak) corresponds to the "Beliandro" on the Tabula Peutingeriana.
}
 \label{fig:PtujCeljeVillach}
 \end{center}
\end{figure}

There is a mountainous river Bel\'a, with the meaning "white", stemming just below 
%the Slovak iconic peak 
Kriv\'a\v n, 
%see Fig. \ref{fig:bela_fella} left, 
one of the highest peaks of the High Tatra mountains and thus of the whole Carpathians. Bel\'a got her name probably by the white smooth rocks in its river-basin. The similar is true for the northern Italian, the Carnic Alpine river Fella,
% see Fig. \ref{fig:bela_fella} right, 
with the Slovenian name Bela. Fella has its spring close to the Friuli-Carinthia-Slovenian border and concerning the current name it went through the standard change of the Slavic "B" to the German "F" recorded now in the Italian official name, see \cite{Miklosic1} and \cite{Miklosic2}-12. This twin Bel\'a-Bela is another clear example of the same geographic name given to mountain rivers by the common Slavic population living in the Alps and Carpathians. Moreover, very near the Fella (Bela) spring, there is the Carinthian town Villach - Beljak in Slovenian. And reading carefully the Tabula Peutingeriana one can find it also there under the name "Beliandro", see Fig. \ref{fig:Peutingeriana2}. We see again the Slavic name of the settlement recorded on the map describing the Roman road system at the beginning of the first millennium AD. We can simply check our claim that Villach (Beljak) corresponds to the "Beliandro" by computing distance from Ptuj ("Petauione" on Tabula Peutingeriana) through Celje ("Celeia") to Villach ("Beliandro"). The whole distance is 137 Roman miles. In order to check the relation of the length of the Roman mile and the Google map distances in kilometres in this area of Tabula Peutingeriana we do it for two clearly given towns and known distances both in kilometres and in the Roman miles. From Ptuj to Celje we have 58,5 km on Google maps and 36 miles on Tabula Peutingeriana, which gives approximate correspondence 1.62 km for 1 Roman mile. Now, 137 x 1.62 km = 222,24 km and its there, see Fig. \ref{fig:PtujCeljeVillach}, the distance from Ptuj to Villach is 220 km on Google maps which perfectly approximate the Tabula Peutingeriana distance 137 Roman miles. Of course, one can check the Villach - "Beliandro" correspondence also in different ways, e.g. by estimating forward distances from "Beliandro" to other stations on the Roman road, but that we left for a reader.

After finding the Bel\'a-Bela twin in the High Tatras and the Carnic Alps one can try to find a corresponding name also in Romania. 
%It was not so straightforward, but we succeeded and it was a nice surprise. 
In the Caraș-Severin county, there is a river 
%with even more interesting Slavic name 
Belareka meaning "Biela rieka" in Slovak and "White river" in English but written with both words together. It is again a mountainous river which joins the river Cerna ("\v Cierna" in Slovak, \cite{Miklosic2}-71, "Black" in English) near the spa town Baile Herculane.
%, see the stylized picture of the town from the year 1824 in Fig. \ref{fig:Belareka}. 
The Cerna river then flows to the Ister and their confluence in nowadays Orșova got the name "Tierua" on Tabula Peutingeriana, see Fig. \ref{fig:Peutingeriana1}. It represents the place of crossing the Ister (lower Danube) by the second Trajan's road to the Dacia province. Interestingly, the place name was also recorded by Klaudios Ptolemaios in Greek as "Δίερνα", see \cite{Safarik, Tierna}. Since in the Greek alphabet there is no direct representation of the letter "\v c", Ptolemaios replaced it by "Δ" due to a pronunciation similarity of \v Cierna-Δίερνα. This record had to be written between AD 85-165, during the Ptolemaios life and after approximately AD 100 when Trajan's roads to Dacia were constructed. We also know by Pavol Jozef \v Saf\'arik \cite{Safarik} that there exists inscription on marble in Mehadia, former Roman camp on Belareka, from AD 157, "Valerius Felix miles coh. IV. stationis TSIERNEN", where "Tsiernen" represents a very accurate transcript of the river Cerna (\v Cierna) name to the Latin.

All these Slavic names writings are dated a few decades later than Strabo's {\it Geographica} was written. Since by Strabo we know that Suevi inhabited that region in the beginning of the first millennium AD, we clearly see the correspondence between the Strabo's Suevi in the Hercynian forest and the Slavs or Slavic settlements in the Carpathian-Alpine basin. We may also claim that the name Suevi (with the Latin "v") gave birth to two nowadays ethnic names of Slavic nations from the former Hercynian Forest, Slovenes and Slovaks and the name form Soebi (with the Greek "b") to the third one, the Serbians which are called "Srbi" in Slovak (in pronunciation quite similar to "Soebi"). And in fact, in the Middle Ages Latin sources up to the end of the 18th century, the Slovak people were called Slavi, Sclavi or gentis Slavae (in pronunciation again very similar to Suevi or Suavi, when considering the term "Suavia" used in \cite{Steinhubel}, page 33), see e.g. "Privilegium pro Slavis solnensis" by \v Ludov\'it I. Ve\v lk\'y (Louis the Great from Anjou) from 1381 where Slovaks are called Sclavi  \cite{Marsina1,Marsina2} or "Historia gentis Slavae" by Juraj (Georgius) Pap\'anek from 1780 \cite{Papanek}.

\begin{figure}[ht]
 \begin{center}
 \includegraphics[width = 1.\textwidth]{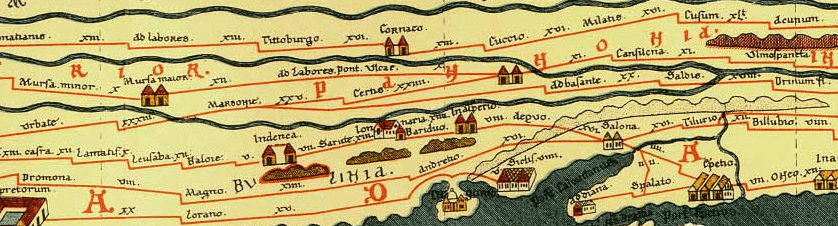}
 \caption{A detail of Tabula Peutingeriana with stations "Vrbate" (left), "pont Vlcae" (middle), "Cerne" (middle) and "Drinum fl." (right) in Bosnia and Slavonia. Source: Wikipedia.
}
 \label{fig:Peutingeriana3}
 \end{center}
\end{figure}

Almost finally we mention "Tergeste" geographic name used by Strabo in {\it Geographica} and appearing also on Tabula Peutingeriana, see Fig. \ref{fig:Peutingeriana2} bottom left. The name "Tergeste" has correspondences in Romania, such as Târgovişte in Dâmbovița County, Muntenia historical region of Romania, which was also the capital of Wallachia between the early 15th century and the 16th century. And the exact correspondence exists also in Slovakia, the village Trhovi\v ste. The common meaning of all these places is the "marketplace" and the name is derived from the Slavic equivalent "trh", "trg" or "targ", see \cite{Miklosic2}-694. There are several further examples such as Târgu Jiu, Târgu Neamț, Târgu Mureș, Târgu Frumos or Târgu Secuiesc in Romania. By the current name similarity, the "Tergeste" used to be related to nowadays Trieste but with some probability, it may also correspond to the city of Monfalcone which is called Tr\v zi\v c in Slovenian with the marketplace meaning \cite{Miklosic2}-694. Concerning Tabula Peutingeriana, we can find several further names of clear Slavic origin in Slavonia and Bosnia regions, see Figure \ref{fig:Peutingeriana3}. There is the station "Vrbate" on the crossing of the Roman road with the nowadays river Vrbas, with the Slavic base of the name, "v\'rba" (in Slovak) meaning the "willow" tree, see \cite{Miklosic2}-746. Then there is the "pont Vlcae" station close to nowadays town Vukovar on the Vuka river, with the Slavic base "vlk" (in Slovak) or "vuk" (in Croatian and Serbian) meaning the "wolf", see \cite{Miklosic2}-733. A further example is the station "Drinum fl" at the place of nowadays Drina river, with the Slavic base "drie\v n" (in Slovak) meaning the "bunchberry", see \cite{Miklosic2}-87. 
%The interesting case is the station "Cerne" which can be localized in the place of nowadays Croatian village Cerna with the meaning "Black" \cite{Miklosic2}-71. However, although its letters are seen quite clearly on the Wikipedia image (Figure \ref{fig:Peutingeriana3}), on the site https://www.tabula-peutingeriana.de/ it is rewritten as "Certis" - why, it is unclear to us. 
The last but not least we mention the karst river Timavo which re-emerges near Monfalcone and empties into the Adriatic. Strabo calls it Timavus in section (5-1-8) of {\it Geographica} and he says it has seven sources by which it re-emerges after 130 stadia in the underground. The Slovenian name of the river is Timava, with the standard Slavic suffix -ava, see \cite{Miklosic1} - Chapter 2, Section V, point 37. And the meaning can be "Tmav\'a" (in Slovak) translated as "dark". This name meaning nicely corresponds to the re-emergence of Timava from the "dark underground". This Slavic name of the river appears also on the Tabula Peutingeriana as "fonte timaui", see Fig. \ref{fig:Peutingeriana2} bottom left, and it again represents a clear example of the Slavic presence in the south-eastern Alpine region at the beginning of the first millennium AD. 

\section{Conclusions} All in all, we have shown the correspondence of the {\it Strabo's Suevi of the Hercynian Forest} with the Slavs in the Carpathian-Alpine basin. In other words, we confirmed the presence of the compact Slavic settlement in this region already at the beginning of the first millennium AD. Interestingly, such a conclusion comes from our mathematical results on the transformation of the Strabo's river Ister to the current map of the world. To perform the transformation we developed a new method for the map registration combining the locally optimal affine transformations with the interpolation/extrapolation by solving numerically the Laplace equation with zero Neumann boundary conditions. The new method keeps the optimality of the local transformations and smoothly interpolate/extrapolate the transformation parameters to the whole computational domain. The method was applied to Strabo's river Ister reconstruction and yield interesting historical conclusions.

The conclusions are in accordance with the opinions and claims of historians and linguists such as Pavol Jozef \v Saf\'arik \cite{Safarik}, \v Ludov\'it \v St\'ur \cite{Stur1,Stur2,Czambel}, Oleg Truba\v cev \cite{Trubacev} or Mario Alinei \cite{Alinei, Alinei2} which all declared the ancientness of the Slavs on the middle Danube, i.e. in the Carpathian-Alpine basin. It supports in some aspects also the narratives of Primary chronicle by the Saint Nestor the Chronicler \cite{Nestor} seeing the middle Danube even as the homeland of all the Slavs, which is however for us too strong and hardly acceptable assumption. More likely, we are in accordance with several aspects of the Paleolithic Continuity Paradigm by Mario Alinei \cite{Alinei, Alinei2} declaring stability of the European population and its ethnic distribution from the very ancient (Upper Paleolithic or Mesolithic) times. Such claims are also in accordance with the results of \cite{Richards} where the authors concluded that the large majority of extant human mitochondrial DNA (mtDNA) lineages entered Europe in several waves during the Late Upper Palaeolithic, mainly around the Last Glacial Maximum, 20000 years ago. After that, the bearers of the mtDNA - sedentary population - adopted in general to a new way of life and the Neolithic demic-diffusion, see e.g. \cite{Renfrew}, and further immigration waves seem to contribute only by about 10-20\% of mtDNA lineages. Also, recent work \cite{Csakyova} of Slovak and Hungarian experts in archaeology and archaeogenetics has shown a large majority of mtDNA haplogroups in medieval (9th-12th century AD) population around Nitra (Western Slovakia) belonging to the Late Upper Palaeolithic period of migration to Europe and similarity of medieval "Slovak" population with nowadays inhabitants of the Carpathian-Alpine basin such as Croats and Romanians. 
%Although the results of this work are based on 12 medieval individuals only, which is statistically rather a small sample, it gives an interesting first view on the subject.
%When reading Strabo carefully and taking into account our result for the Ister course, we see such ethnic stability in Central Europe quite clearly, of course with some exceptions caused mainly by a centralized education in the last centuries. 

We hope that results of this paper can help geographers and cartographers with the historical map registration and 
%at least slightly 
also to historians
%, such as J\'an Steinh\"ubel, 
to localize the "Suavia" at the beginning of the first millennium AD \cite{Steinhubel} or 
%to Florin Curta \cite{Curta} and Martin Homza \cite{Homza} 
to explain how it was possible that Slavs appeared so fast in the Carpathian basin and north of the Istros (Lower Danube) at the beginning of the second half of the first millennium AD \cite{Curta,Homza}. In fact, the Slavs did not appear, they were there "since time immemorial".

%\section{REFERENCES}
%
%
%\bibliographystyle{misas}%%%{plain}%%{alpha}%{misas} ...
%\bibliography{tmmp000}

%
\end{document}